\documentclass[12pt]{amsart}
\usepackage[colorlinks=true, pdfstartview=FitV, linkcolor=blue, citecolor=blue, urlcolor=blue, breaklinks=true]{hyperref}
\usepackage{amsmath,amsfonts,amssymb,amsthm,amscd,comment}
\usepackage{tikz}
\usepackage{mathtools}
\usetikzlibrary{decorations.markings}

\DeclareMathOperator{\dimbar}{\underbar{$\mathrm{dim}$}}
\DeclareMathOperator{\Hom}{Hom}
\DeclareMathOperator{\Ext}{Ext^1}
\DeclareMathOperator{\Aut}{Aut}

\DeclareMathOperator{\vspan}{span}
\DeclareMathOperator{\ch}{char}
\DeclareMathOperator{\id}{id}
\DeclareMathOperator{\supp}{supp}
\DeclareMathOperator{\diag}{diag}
\DeclareMathOperator{\rad}{rad}
\DeclareMathOperator{\gldim}{gl.dim}
\DeclareMathOperator{\avq}{\mathfrak{Q}_\mathrm{abs}}
\DeclareMathOperator{\rvq}{\mathfrak{Q}_\mathrm{rel}}
\DeclareMathOperator{\spec}{\mathfrak{M}}

\DeclareMathOperator{\ad}{ad}
\DeclareMathOperator{\f}{\mathbf{f}}
\DeclareMathOperator{\g}{\mathfrak{g}}

\newcommand\Q{\mathcal{Q}}
\newcommand\U{\mathbf{U}}
\newcommand\F{\mathbf{F}}
\newcommand\R{\mathfrak{R}}
\newcommand\HH{\mathcal{H}}
\newcommand\A{\mathcal{A}}

\theoremstyle{plain}
\newtheorem{theo}{Theorem}[section]
\newtheorem{prop}[theo]{Proposition}
\newtheorem{lem}[theo]{Lemma}
\newtheorem{cor}[theo]{Corollary}

\newtheorem{conj}[theo]{Conjecture}

\theoremstyle{definition}
\newtheorem{defn}[theo]{Definition}
\newtheorem{rmk}[theo]{Remark}
\newtheorem{egg}[theo]{Example}

\numberwithin{equation}{section}

\allowdisplaybreaks


\begin{document}

\title[Valued Graphs and the Representation Theory of Lie Algebras]{Valued Graphs and the Representation Theory \linebreak of Lie Algebras}

\author{Joel Lemay}
\address{J.~Lemay: Department of Mathematics and Statistics, University of Ottawa}


\begin{abstract}
 Quivers (directed graphs) and species (a generalization of quivers) and their representations play a key role in many areas of mathematics including combinatorics, geometry, and algebra. Their importance is especially apparent in their applications to the representation theory of associative algebras, Lie algebras, and quantum groups. In this paper, we discuss the most important results in the representation theory of species, such as Dlab and Ringel's extension of Gabriel's theorem, which classifies all species of finite and tame representation type. We also explain the link between species and $K$-species (where $K$ is a field). Namely, we show that the category of $K$-species can be viewed as a subcategory of the category of species. Furthermore, we prove two results about the structure of the tensor ring of a species containing no oriented cycles that do not appear in the literature. Specifically, we prove that two such species have isomorphic tensor rings if and only if they are isomorphic as ``crushed'' species, and we show that if $K$ is a perfect field, then the tensor algebra of a $K$-species tensored with the algebraic closure of $K$ is isomorphic to, or Morita equivalent to, the path algebra of a quiver.
\end{abstract}

\maketitle

\tableofcontents

\section*{Introduction}
Species and their representations were first introduced in 1973 by Gabriel in \cite{g}. Let $K$ be a field. Let $A$ be a finite-dimensional, associative, unital, basic $K$-algebra and let $\rad A$ denote its Jacobson radical. Then $A/ \rad A \cong \Pi_{i \in \mathcal{I}} K_i$, where $\mathcal{I}$ is a finite set and $K_i$ is a finite-dimensional $K$-division algebra for each $i \in \mathcal{I}$. Moreover, $\rad A / (\rad A)^2 \cong \bigoplus_{i,j \in \mathcal{I}} {_jM}_i$, where ${_jM}_i$ is a finite-dimensional $(K_j,K_i)$-bimodule for each $i, j \in \mathcal{I}$. We then associate to $A$ a valued graph $\Delta_A$ with vertex set $\mathcal{I}$ and valued arrows $i \xrightarrow{(d_{ij},d_{ji})} j$ for each ${_jM}_i \ne 0$, where $d_{ij}=\dim_{K_j}({_jM}_i)$. The valued graph $\Delta_A$, the division algebras $K_i$ ($i \in \mathcal{I}$) and the bimodules $_jM_i$ ($i, j \in \mathcal{I}$) constitute a \emph{species} and contain a great deal of information about the representation theory of $A$ (in some cases, \emph{all} the information). When working over an algebraically closed field, a species is simply a quiver (directed graph) in the sense that all $K_i \cong K$ and all ${_jM}_i \cong K^n$ so only $\Delta_A$ is significant. In this case, Gabriel was able to classify all quivers of finite representation type (that is, quivers with only finitely many non-isomorphic indecomposable representations); they are precisely those whose underlying graph is a (disjoint union of) Dynkin diagram(s) of type A, D or E. Moreover, he discovered that the isomorphism classes of indecomposable representations of these quivers are in bijection with the positive roots of the Kac-Moody Lie algebra associated to the corresponding diagram. Gabriel's theorem is the starting point of a series of remarkable results such as the construction of Kac-Moody Lie algebras and quantum groups via Ringel-Hall algebras, the geometry of quiver varieties, and Lusztig's categorification of quantum groups via perverse sheaves. Lusztig, for example, was able to give a geometric interpretation of the positive part of quantized enveloping algebras using quiver varieties (see \cite{lusz}).

While quivers are useful tools in representation theory, they have their limitations. In particular, their application to the representation theory of associative unital algebras, in general, only holds when working over an algebraically closed field. Moreover, the Lie theory that is studied by quiver theoretic methods is naturally that of symmetric Kac-Moody Lie algebras. However, many of the fundamental examples of Lie algebras of interest to mathematicians and physicists are symmetrizable Kac-Moody Lie algebras which are not symmetric. Species allow us to relax these limitations.

In his paper \cite{g}, Gabriel outlined how one could classify all species of finite representation type over non-algebraically closed fields. However, it was Dlab and Ringel in 1976 (see \cite{dr}) who were ultimately able to generalize Gabriel's theorem and show that a species is of finite representation type if and only if its underlying valued graph is a Dynkin diagram of finite type. They also showed that, just as for quivers, there is a bijection between the isomorphism classes of the indecomposable representations and the positive roots of the corresponding Kac-Moody Lie algebra.

Despite having been introduced at the same time, the representation theory of quivers is much more well-known and well-developed than that of species. In fact, the very definition of species varies from text to text; some use the more ``general'' definition of a species (e.g.\ \cite{dr}) while others use the alternate definition of a $K$-species (e.g.\ \cite{ddpw}). Yet the relationship between these two definitions is rarely discussed. Moreover, while there are many well-known results in the representation theory of quivers, such as Gabriel's theorem or Kac's theorem, it is rarely mentioned whether or not these results generalize for species. Indeed, there does not appear to be any single comprehensive reference for species in the literature. The main goal of this paper is to compare the current literature and collect all the major, often hard to find, results in the representation theory of species into one text.

This paper is divided into seven sections. In the first, we give all the preliminary material on quivers and valued quivers that will be needed for the subsequent sections. In particular, we address the fact that two definitions of valued quivers exist in the literature.

In Section 2, we define both species and $K$-species and discuss how the definitions are related. Namely, we define the categories of species and $K$-species and show that the category of $K$-species can be thought of as a subcategory of the category of species. That is, via an appropriate functor, all $K$-species are species. There are, however, species that are not $K$-species for any field $K$.

The third section deals with the tensor ring (resp.\ algebra) $T(\Q)$ associated to a species (resp.\ $K$-species) $\Q$. This is a generalization of the path algebra of a quiver. If $K$ is a perfect field, then for any finite-dimensional associative unital $K$-algebra $A$, the category of $A$-modules is equivalent to the category of $T(\Q)/I$-modules for some $K$-species $\Q$ and some ideal $I$. Also, it will be shown in Section 6, that the category of representations of $\Q$ is equivalent to the category of $T(\Q)$-modules. These results show why species are such important tools in representation theory; modulo an ideal, they allow us to understand the representation theory of finite-dimensional associative unital algebras.

In Section 4, we follow the work of \cite{ddpw} to show that, when working over a finite field, one can simply deal with quivers (with automorphism) rather than species. That is, we show that if $\Q$ is an $\mathbb{F}_q$-species, then the tensor algebra of $\Q$ is isomorphic to the fixed point algebra of the path algebra of a quiver under the Frobenius morphism.

In the fifth section, we further discuss the link between a species and its tensor ring. In particular, we prove two results that do not seem to appear in the literature.

\mbox{}

\begin{raggedleft}\textbf{Theorem \ref{theo:iso}}\end{raggedleft} \textit{Let $\Q$ and $\Q'$ be two species with no oriented cycles. Then $T(\Q) \cong T(\Q')$ if and only if $\Q^C \cong \Q'^C$ (where $\Q^C$ and $\Q'^C$ denote the \emph{crushed} species of $\Q$ and $\Q'$).}

\mbox{}

\begin{raggedleft}\textbf{Theorem \ref{theo:tensor} and Corollary \ref{cor:tensor}}\end{raggedleft} \textit{Let $K$ be a perfect field and $\Q$ a $K$-species containing no oriented cycles. Then $\overline{K} \otimes_K T(\Q)$ is isomorphic to, or Morita equivalent to, the path algebra of a quiver (where $\overline{K}$ denotes the algebraic closure of $K$).}

\mbox{}

Section 6 deals with representations of species. We discuss many of the most important results in the representation theory of quivers, such as the theorems of Gabriel and Kac, and their generalizations for species.

The seventh and final section deals with the Ringel-Hall algebra of a species. It is well-known that the generic composition algebra of a quiver is isomorphic to the positive part of the quantized enveloping algebra of the associated Kac-Moody Lie algebra. Also, Sevenhant and Van Den Bergh have shown that the Ringel-Hall algebra itself is isomorphic to the positive part of the quantized enveloping algebra of a generalized Kac-Moody Lie algebra (see \cite{sv}). We show that these results hold for species as well. While this is not a new result, it does not appear to be explained in detail in the literature.

We assume throughout that all algebras (other than Lie algebras) are associative and unital.

\medskip

\paragraph{\textbf{Acknowledgements}} First and foremost, I would like to thank Prof.\ Alistair Savage for introducing me to this topic and for his invaluable guidance and encouragement. Furthermore, I would like to thank Prof.\ Erhard Neher and Prof.\ Vlastimil Dlab for their helpful comments and advice. 

\section{Valued Quivers}
In this section, we present the preliminary material on quivers and valued quivers that will be used throughout this paper. In particular, we begin with the definition of a quiver and then discuss valued quivers. There are two definitions of valued quivers that can be found in the literature; we present both and give a precise relationship between the two in terms of a functor between categories (see Lemma \ref{lem:functorF}). We also discuss the idea of ``folding'', which allows one to obtain a valued quiver from a quiver with automorphism.

\begin{defn}[Quiver] 
  A \emph{quiver}\index{Quiver} $Q$ is a directed graph. That is, $Q=(Q_0,Q_1,t,h)$, where $Q_0$ and $Q_1$ are sets and $t$ and $h$ are set maps $Q_1 \rightarrow Q_0$. The elements of $Q_0$ are called \emph{vertices} and the elements of $Q_1$ are called \emph{arrows}. For every $\rho \in Q_1$, we call $t(\rho)$ the \emph{tail} of $\rho$ and $h(\rho)$ the \emph{head} of $\rho$. By an abuse of notation, we often simply write $Q=(Q_0,Q_1)$ leaving the maps $t$ and $h$ implied. The sets $Q_0$ and $Q_1$ may well be infinite; however we will deal exclusively with quivers having only finitely many vertices and arrows. We will also restrict ourselves to quivers whose underlying undirected graphs are connected.

  A \emph{quiver morphism} $\varphi: Q \rightarrow Q'$ consists of two set maps, $\varphi_0: Q_0 \rightarrow Q'_0$ and $\varphi_1: Q_1 \rightarrow Q'_1$, such that $\varphi_0(t(\rho))=t(\varphi_1(\rho))$ and $\varphi_0(h(\rho))=h(\varphi_1(\rho))$ for each $\rho \in Q_1$.
\end{defn}

For $\rho \in Q_1$, we will often use the notation $\rho:i \rightarrow j$ to mean $t(\rho)=i$ and $h(\rho)=j$.

\begin{defn}[Absolute valued quiver] 
  An \emph{absolute valued quiver}\index{Quiver!valued!absolute} is a quiver $\Gamma=(\Gamma_0,\Gamma_1)$ along with a positive integer $d_i$ for each $i \in \Gamma_0$ and a positive integer $m_\rho$ for each $\rho \in \Gamma_1$ such that $m_\rho$ is a common multiple of $d_{t(\rho)}$ and $d_{h(\rho)}$ for each $\rho \in \Gamma_1$. We call $(d_i,m_\rho)_{i \in \Gamma_0, \rho \in \Gamma_1}$ an \emph{(absolute) valuation}\index{Valuation!absolute} of $\Gamma$. By a slight abuse of notation, we often refer to $\Gamma$ as an absolute valued quiver, leaving the valuation implied.

  An \emph{absolute valued quiver morphism} is a quiver morphism $\varphi: \Gamma \rightarrow \Gamma'$ respecting the valuations. That is, $d'_{\varphi_0(i)}=d_i$ for each $i \in \Gamma_0$ and $m'_{\varphi_1(\rho)}=m_\rho$ for each $\rho \in \Gamma_1$.

  Let $\avq$\index{$\avq$} denote the category of absolute valued quivers.
\end{defn}

A (non-valued) quiver can be viewed as an absolute valued quiver with trivial values (i.e.\ all $d_i=m_\rho=1$). Thus, valued quivers are a generalization of quivers.

Given a quiver $Q$ and an automorphism $\sigma$ of $Q$, we can construct an absolute valued quiver $\Gamma$ with valuation $(d_i,m_\rho)_{i\in \Gamma_0, \rho \in \Gamma_1}$ by ``folding''\index{Folding} $Q$ as follows:

\begin{itemize}
  \item $\Gamma_0=\{\text{vertex orbits of } \sigma\}$,
  \item $\Gamma_1=\{\text{arrow orbits of } \sigma\}$,
  \item for each $i \in \Gamma_0$, $d_i$ is the number of vertices in the orbit $i$,
  \item for each $\rho \in \Gamma_1$, $m_\rho$ is the number of arrows in the orbit $\rho$.
\end{itemize}

Given $\rho \in \Gamma_1$, let $m=m_\rho$ and $d=d_{t(\rho)}$. The orbit $\rho$ consists of $m$ arrows in $Q_0$, say $\{\rho_i=\sigma^{i-1}(\rho_1)\}_{i=1}^m$. Because $\sigma$ is a quiver morphism, we have that each $t(\rho_i)=t(\sigma^{i-1}(\rho_1))$ is in the orbit $t(\rho)$ and that $t(\rho_i)=t(\sigma^{i-1}(\rho_1))=\sigma^{i-1}(t(\rho_1))$. The value $d$ is the least positive integer such that $\sigma^d(t(\rho_1))=t(\rho_1)$ and since $\sigma^m(t(\rho_1))=t(\rho_1)$ (because $\sigma^m(\rho_1)=\rho_1$), then $d\;|\;m$. By the same argument, $d_{h(\rho)}\;|\;m$. Thus, this construction does in fact yield an absolute valued quiver.

Conversely, given an absolute valued quiver $\Gamma$ with valuation $(d_i,m_\rho)_{i \in \Gamma_0,\rho \in \Gamma_1}$, it is possible to construct a quiver with automorphism $(Q,\sigma)$ that folds into $\Gamma$ in the following way. Let $\overline{x}_y$ be the unique representative of $(x$ mod $y)$ in the set $\{1, 2, \dots, y\}$ for $x,y$ positive integers. Then define:

\begin{itemize}
  \item $Q_0=\{v_i(j) \; | \; i \in \Gamma_0, \; 1 \le j \le d_i \}$,

  \item $Q_1=\{a_\rho(k) \; | \; \rho \in \Gamma_1, \; 1 \le k \le m_\rho \}$,

  \item $t \left(a_\rho(k) \right) = v_{t(\rho)}\left(\overline{k}_{d_{t(\rho)}} \right) $ and $h \left(a_\rho(k) \right) = v_{h(\rho)}\left(\overline{k}_{d_{h(\rho)}} \right)$,

  \item $\sigma \left(v_i(j)\right) = v_i\left( \overline{(j+1)}_{d_i} \right)$,

  \item $\sigma \left(a_\rho(k)\right) = a_\rho\left( \overline{(k+1)}_{m_\rho} \right)$.
\end{itemize}
It is clear that $Q$ is a quiver. It is easily verified that $\sigma$ is indeed an automorphism of $Q$. Given the construction, we see that $(Q,\sigma)$ folds into $\Gamma$. However, we do not have a one-to-one correspondence between absolute valued quivers and quivers with automorphism since, in general, several non-isomorphic quivers with automorphism can fold into the same absolute valued quiver, as the following example demonstrates.

\begin{egg}\cite[Example 3.4]{ddpw} Consider the following two quivers.

  \begin{center}

    \begin{tikzpicture} [>=stealth]
      \draw (-1,0) node {$Q:$} (4,0) node {$Q':$};

      \filldraw (0,0) circle (2pt)	(5,0) circle (2pt)
		(1,0.7) circle (2pt)	(6,0.7) circle (2pt)
		(1,-0.7) circle (2pt)	(6,-0.7) circle (2pt)
		(2,0.7) circle (2pt)	(7,0.7) circle (2pt)
		(2,-0.7) circle (2pt)	(7,-0.7) circle (2pt);

      \draw [->] (0,0) -- (0.6,0.42); \draw (0.6,0.42) -- (1,0.7);
      \draw [->] (0,0) -- (0.6,-0.42); \draw (0.6,-0.42) -- (1,-0.7);
      \draw [->] (1,0.75) -- (1.6,0.75); \draw (1.6,0.75) -- (2,0.75);
      \draw [->] (1,0.65) -- (1.6,0.65); \draw (1.6,0.65) -- (2,0.65);
      \draw [->] (1,-0.75) -- (1.6,-0.75); \draw (1.6,-0.75) -- (2,-0.75);
      \draw [->] (1,-0.65) -- (1.6,-0.65); \draw (1.6,-0.65) -- (2,-0.65);

      \draw [->] (5,0) -- (5.6,0.42); \draw (5.6,0.42) -- (6,0.7);
      \draw [->] (5,0) -- (5.6,-0.42); \draw (5.6,-0.42) -- (6,-0.7);
      \draw [->] (6,0.7) -- (6.6,0.7); \draw (6.6,0.7) -- (7,0.7);
      \draw [->] (6,0.7) -- (6.6,-0.14); \draw (6.6,-0.14) -- (7,-0.7);
      \draw [->] (6,-0.7) -- (6.6,-0.7); \draw (6.6,-0.7) -- (7,-0.7);
      \draw [->] (6,-0.7) -- (6.6,0.14); \draw (6.6,0.14) -- (7,0.7);

      \begin{small}
	\draw (-0.3,0) node {1} (1,1) node {2} (2,1) node {3} (1,-1) node {4} (2,-1) node {5};

	\draw (4.7,0) node {$a$} (6,1) node {$b$} (7,1) node {$c$} (6,-1) node {$d$} (7,-1) node {$e$};

	\draw (0.5,0.6) node {$\alpha_1$} (0.5,-0.6) node {$\alpha_2$} (1.5,0.95) node {$\alpha_3$} (1.5,0.45) node {$\alpha_4$} (1.5,-0.45) node {$\alpha_5$} (1.5,-0.95) node {$\alpha_6$};

	\draw (5.5,0.6) node {$\beta_1$} (5.5,-0.65) node {$\beta_2$} (6.5,0.9) node {$\beta_3$} (7,-0.3) node {$\beta_4$} (7,0.3) node {$\beta_5$} (6.5,-0.95) node {$\beta_6$};
      \end{small}
    \end{tikzpicture}

  \end{center}
  Define $\sigma \in \Aut(Q)$ and $\sigma' \in \Aut(Q')$ by
  \[\sigma:   \begin{pmatrix}
	      1 & 2 & 3 & 4 & 5 \\
	      1 & 4 & 5 & 2 & 3
             \end{pmatrix}, \;
	     \begin{pmatrix}
	      \alpha_1 & \alpha_2 & \alpha_3 & \alpha_4 & \alpha_5 & \alpha_6 \\
	      \alpha_2 & \alpha_1 & \alpha_5 & \alpha_6 & \alpha_3 & \alpha_4
	     \end{pmatrix} \]

  \[\sigma':  \begin{pmatrix}
	      a & b & c & d & e \\
	      a & d & e & b & c
             \end{pmatrix}, \;
	     \begin{pmatrix}
	      \beta_1 & \beta_2 & \beta_3 & \beta_4 & \beta_5 & \beta_6 \\
	      \beta_2 & \beta_1 & \beta_6 & \beta_5 & \beta_4 & \beta_3
	     \end{pmatrix}.\]
  Then, both $(Q,\sigma)$ and $(Q',\sigma')$ fold into

  
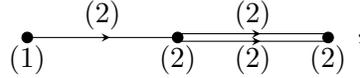
\begin{center}

    \begin{tikzpicture} [>=stealth, baseline=0]
      \filldraw (0,0) circle (2pt) (2,0) circle (2pt) (4,0) circle (2pt);

      \draw [->] (0,0) -- (1.1,0); \draw (1,0) -- (2,0);
      \draw [->] (2,0.05) -- (3.1,0.05); \draw (3,0.05) -- (4,0.05);
      \draw [->] (2,-0.05) -- (3.1,-0.05); \draw (3,-0.05) -- (4,-0.05);

      \begin{small}
	\draw (0,-0.3) node {(1)} (2,-0.3) node {(2)} (4,-0.3) node {(2)};
	\draw (1,0.3) node {(2)} (3,0.3) node {(2)} (3,-0.3) node {(2)};
      \end{small}
    \end{tikzpicture},

  \end{center}
  yet $Q$ and $Q'$ are not isomorphic as quivers. 

\end{egg}

\begin{defn}[Relative valued quiver]\label{defn:rvq} 
  A \emph{relative valued quiver}\index{Quiver!valued!relative} is a quiver $\Delta=(\Delta_0, \Delta_1)$ along with positive integers $d^\rho_{ij}$, $d^\rho_{ji}$ for each arrow $\rho:i \rightarrow j$ in $\Delta_1$ such that there exist positive integers $f_i$, $i \in \Delta_0$, satisfying
  \[
    d^\rho_{ij}f_j = d^\rho_{ji}f_i
  \]
  for all arrows $\rho:i \rightarrow j$ in $\Delta_1$. We call $(d^\rho_{ij},d^\rho_{ji})_{(\rho: i \rightarrow j) \in \Delta_1}$ a \emph{(relative) valuation}\index{Valuation!relative} of $\Delta$. By a slight abuse of notation, we often refer to $\Delta$ as a relative valued quiver, leaving the valuation implied.

  We will use the notation:

  \begin{center}

    \begin{tikzpicture} [>=stealth, baseline=0]

      \filldraw (0,0) circle (2pt) (1.5,0) circle (2pt);

      \draw [->] (0,0) -- (0.8,0); \draw (0.7,0) -- (1.5,0);

      \begin{small}
	\draw (0,-0.3) node {$i$} (1.5,-0.3) node {$j$} (0.75,-0.3) node {$\rho$} (0.75,0.3) node {$(d^\rho_{ij},d^\rho_{ji})$};
      \end{small}
    \end{tikzpicture}.

  \end{center}
  In the case that $(d^\rho_{ij},d^\rho_{ji})=(1,1)$, we simply omit it.

  A \emph{relative valued quiver morphism} is a quiver morphism $\varphi: \Delta \rightarrow \Delta'$ satisfying:
  \[ (d')^{\varphi_1(\rho)}_{\varphi_0(i)\varphi_0(j)} = d^\rho_{ij} \quad \text{ and } \quad (d')^{\varphi_1(\rho)}_{\varphi_0(j)\varphi_0(i)} = d^\rho_{ji}\]
  for all arrows $\rho:i \rightarrow j$ in $\Delta_1$.

  Let $\rvq$\index{$\rvq$} denote the category of relative valued quivers.

\end{defn}

Note that the definition of a relative valued quiver closely resembles the definition of a symmetrizable Cartan matrix. We will explore the link between the two in Section 6, which deals with representations.

As with absolute valued quivers, one can view (non-valued) quivers as relative valued quivers with trivial values (i.e.\ all $(d^\rho_{ij},d^\rho_{ji})=(1,1)$). Thus, relative valued quivers are also a generalization of quivers.

It is natural to ask, then, how the two categories $\avq$ and $\rvq$ are related. Given $\Gamma \in \avq$ with valuation $(d_i,m_\rho)_{i \in \Gamma_0, \rho \in \Gamma_1}$, define $\F(\Gamma) \in \rvq$\index{$\F$} with valuation $(d^\rho_{ij},d^\rho_{ji})_{(\rho: i \rightarrow j) \in \F(\Gamma)_1}$ as follows:

\begin{itemize}

  \item the underlying quiver of $\F(\Gamma)$ is equal to that of $\Gamma$,
  \item the values $(d^\rho_{ij},d^\rho_{ji})$ are given by:
    \[ d^\rho_{ij}=\dfrac{m_\rho}{d_j} \quad \text{ and } \quad d^\rho_{ji}=\dfrac{m_\rho}{d_i} \]
  for all arrows $\rho:i \rightarrow j$ in $\F(\Gamma)_1$.

\end{itemize}

It is clear that $\F(\Gamma)$ satisfies the definition of a relative valued quiver (simply set all the $f_i=d_i$). Given a morphism $\varphi: \Gamma \rightarrow \Gamma'$ in $\avq$, one can simply define $\F(\varphi): \F(\Gamma) \rightarrow \F(\Gamma')$ to be the morphism given by $\varphi$, since $\Gamma$ and $\Gamma'$ have the same underlying quivers as $\F(\Gamma)$ and $\F(\Gamma')$, respectively. By construction of $\F(\Gamma)$ and $\F(\Gamma')$, it is clear then that $\F(\varphi)$ is a morphism in $\rvq$. Thus, $\F$ is a functor from $\avq$ to $\rvq$.

\begin{lem}\label{lem:functorF}

  The functor $\F: \avq \rightarrow \rvq$ is faithful and surjective.

\end{lem}

\begin{proof}
  Suppose $\F(\varphi) = \F(\psi)$ for two morphisms $\varphi, \psi: \Gamma \rightarrow \Gamma'$ in $\avq$. By definition, $\F(\varphi) = \varphi$ on the underlying quivers of $\Gamma$ and $\Gamma'$. Likewise for $\F(\psi)$ and $\psi$. Thus, $\varphi = \psi$ and $\F$ is faithful.

  Suppose $\Delta$ is a relative valued quiver. By definition, there exist positive integers $f_i$, $i \in \Delta_0$, such that $d^\rho_{ij}f_j = d^\rho_{ji}f_i$ for each arrow $\rho: i \rightarrow j$ in $\Delta_1$. Fix a particular choice of these $f_i$. Define $\Gamma \in \avq$ as follows:

  \begin{itemize}
    \item the underlying quiver of $\Gamma$ is the same as that of $\Delta$,
    \item set $d_i = f_i$ for each $i \in \Gamma_0=\Delta_0$,
    \item set $m_\rho = d^\rho_{ij}f_j = d^\rho_{ji}f_i$ for each arrow $\rho:i \rightarrow j$ in $\Gamma_1 = \Delta_1$.
  \end{itemize}
  Then, $\Gamma$ is an absolute valued quiver and $\F(\Gamma) = \Delta$. Thus, $\F$ is surjective.

\end{proof}

Note that $\F$ is not full, and thus not an equivalence of categories, as the following example illustrates.

\begin{egg}

  Consider the following two non-isomorphic absolute valued quivers.

  \begin{center}

    \begin{tikzpicture} [>=stealth]
      \draw (-1,0) node {$\Gamma$:} (3,0) node {$\Gamma'$:};

      \filldraw (0,0) circle (2pt) (1,0) circle (2pt);

      \filldraw (4,0) circle (2pt) (5,0) circle (2pt);

      \draw [->] (0,0) -- (0.6,0); \draw (0.6,0) -- (1,0);

      \draw [->] (4,0) -- (4.6,0); \draw (4.6,0) -- (5,0);

      \begin{small}
	\draw (0,-0.3) node {(2)} (0.5,0.3) node {(2)} (1,-0.3) node {(1)};

	\draw (4,-0.3) node {(4)} (4.5,0.3) node {(4)} (5,-0.3) node {(2)};
      \end{small}
    \end{tikzpicture}

  \end{center}
  Both $\Gamma$ and $\Gamma'$ are mapped to:

    \begin{center}

      \begin{tikzpicture} [>=stealth, baseline=0]
      \draw (-1.7,0) node {$\F(\Gamma)=\F(\Gamma')$:};

      \filldraw (0,0) circle (2pt) (1,0) circle (2pt);

      \draw [->] (0,0) -- (0.6,0); \draw (0.6,0) -- (1,0);

      \begin{small}
	\draw (0.5,0.3) node {(2,1)};
      \end{small}
    \end{tikzpicture}$\;$.

    \end{center}
  One sees that $\Hom_{\avq}(\Gamma, \Gamma')$ is empty whereas $\Hom_{\rvq}(\F(\Gamma),\F(\Gamma'))$ is not (it contains the identity). Thus, \[\F:\Hom_{\avq}(\Gamma, \Gamma') \rightarrow \Hom_{\rvq}(\F(\Gamma),\F(\Gamma'))\] is not surjective, and hence $\F$ is not full.

\end{egg}

It is tempting to think that one could remedy this by restricting $\F$ to the full subcategory of $\avq$ consisting of objects $\Gamma$ with valuations $(d_i,m_\rho)_{i \in \Gamma_0,\rho \in \Gamma_1}$ such that the greatest common divisor of all $d_i$ is 1. While one can show that $\F$ restricted to this subcategory is injective on objects, it would still not be full, as the next example illustrates.

\begin{egg}\label{egg:notfull}

  Consider the following two absolute valued quivers.

  \begin{center}

    \begin{tikzpicture} [>=stealth]
      \draw (-1,0) node {$\Gamma$:} (3,0) node {$\Gamma'$:};

      \filldraw (0,0) circle (2pt) (1,0) circle (2pt);

      \filldraw (4,0) circle (2pt) (5,0) circle (2pt) (6,0) circle (2pt);

      \draw [->] (0,0) -- (0.6,0); \draw (0.6,0) -- (1,0);

      \draw [->] (4,0) -- (4.6,0); \draw (4.6,0) -- (5,0); \draw [->] (5,0) -- (5.6,0); \draw (5.6,0) -- (6,0);

      \begin{small}
	\draw (0,-0.3) node {(2)} (0.5,0.3) node {(2)} (1,-0.3) node {(1)}; \draw (0.5,-0.2) node{$\rho$};

	\draw (4,-0.3) node {(4)} (4.5,0.3) node {(4)} (5,-0.3) node {(2)} (5.5,0.3) node {(2)} (6,-0.3) node {(1)}; \draw (4.5,-0.2) node {$\alpha$} (5.5,-0.3) node {$\beta$};
      \end{small}
    \end{tikzpicture}

  \end{center}
  The values of the vertices of $\Gamma$ have greatest common divisor 1. The same is true of $\Gamma'$. By applying $\F$ we get:

  \begin{center}

    \begin{tikzpicture} [>=stealth]
      \draw (-1,0) node {$\F(\Gamma)$:} (3,0) node {$\F(\Gamma')$:};

      \filldraw (0,0) circle (2pt) (1,0) circle (2pt);

      \filldraw (4,0) circle (2pt) (5,0) circle (2pt) (6,0) circle (2pt);

      \draw [->] (0,0) -- (0.6,0); \draw (0.6,0) -- (1,0);

      \draw [->] (4,0) -- (4.6,0); \draw (4.6,0) -- (5,0); \draw [->] (5,0) -- (5.6,0); \draw (5.6,0) -- (6,0);

      \begin{small}
	\draw (0.5,0.3) node {(2,1)}; \draw (0.5,-0.2) node{$\rho$};

	\draw (4.5,0.3) node {(2,1)} (5.5,0.3) node {(2,1)}; \draw (4.5,-0.2) node {$\alpha$} (5.5,-0.3) node {$\beta$};
      \end{small}
    \end{tikzpicture}

  \end{center}
  One sees that $\Hom_{\avq}(\Gamma, \Gamma')$ contains only one morphism (induced by $\rho \mapsto \beta$), while on the other hand $\Hom_{\rvq}(\F(\Gamma),\F(\Gamma'))$ contains two morphisms (induced by $\rho \mapsto \alpha$ and $\rho \mapsto \beta$). Thus, \[\F:\Hom_{\avq}(\Gamma, \Gamma') \rightarrow \Hom_{\rvq}(\F(\Gamma),\F(\Gamma'))\] is not surjective, and hence $\F$ is not full, even when restricted to the subcategory of objects with vertex values having greatest common divisor 1.
\end{egg}

Note that there is no similar functor $\rvq \rightarrow \avq$. Following the proof of Lemma \ref{lem:functorF}, one sees that finding a preimage under $\F$ of a relative valued quiver $\Delta$ is equivalent to making a choice of $f_i$ (from Definition \ref{defn:rvq}). One can show that there is a unique such choice satisfying $\gcd(f_i)_{i \in \Delta_0} = 1$ (so long as $\Delta$ is connected). Thus, there is a natural and well-defined way to map objects of $\rvq$ to objects of $\avq$ by mapping a relative valued quiver $\Delta$ to the unique absolute valued quiver $\Gamma$ with valuation $(d_i,m_\rho)_{i \in \Gamma_0, \rho \in \Gamma_1}$ satisfying $\F(\Gamma)=\Delta$ and $\gcd(d_i)_{i \in \Gamma_0}=1$. However, there is no such natural mapping on the morphisms of $\rvq$. For instance, under this natural mapping on objects, in Example \ref{egg:notfull}, the relative valued quivers $\F(\Gamma)$ and $\F(\Gamma')$ are mapped to $\Gamma$ and $\Gamma'$, respectively. However, there is no natural way to map the morphism $\F(\Gamma) \rightarrow \F(\Gamma')$ induced by $\rho \mapsto \alpha$ to a morphism $\Gamma \rightarrow \Gamma'$ since there is no morphism $\Gamma \rightarrow \Gamma'$ such that $\rho \mapsto \alpha$. Thus, there does not appear to be a functor similar to $\F$ from $\rvq$ to $\avq$.

\section{Species and $K$-species}
The reason for introducing two different definitions of valued quivers in the previous section, is that there are two different definitions of species in the literature: one for each of the two versions of valued quivers. In this section, we introduce both definitions of species and discuss how they are related (see Proposition \ref{prop:species} as well as Examples \ref{egg:species1}, \ref{egg:species2} and \ref{egg:species3}).

First, we begin with the more general definition of species (see for example \cite{g} or \cite{dr}). Recall that if $R$ and $S$ are rings and $M$ is an $(R,S)$-bimodule, then $\Hom_R(M,R)$ is an $(S,R)$-bimodule via $(s \cdot \varphi \cdot r)(m) = \varphi(m \cdot s) r$ and $\Hom_S(M,S)$ is an $(S,R)$-bimodule via $(s \cdot \psi \cdot r)(m) = s \psi(r \cdot m)$ for all $r \in R$, $s \in S$, $m \in M$, $\varphi \in \Hom_R(M,R)$ and $\psi \in \Hom_S(M,S)$.

\begin{defn}[Species] \label{defn:spec} 
  Let $\Delta$ be a relative valued quiver\index{Quiver!valued!relative} with valuation \linebreak $(d^\rho_{ij},d^\rho_{ji})_{(\rho:i\rightarrow j) \in \Delta_1}$. A \emph{modulation}\index{Modulation} $\mathbb{M}$ of $\Delta$ consists of a division ring $K_i$ for each $i \in \Delta_0$, and a $(K_{h(\rho)},K_{t(\rho)})$-bimodule $M_\rho$ for each $\rho \in \Delta_1$ such that the following two conditions hold:

  \begin{enumerate}
    \item $\Hom_{K_{t(\rho)}}(M_\rho,K_{t(\rho)}) \cong \Hom_{K_{h(\rho)}}(M_\rho,K_{h(\rho)})$ as $(K_{t(\rho)},K_{h(\rho)})$-bimodules, and
    \item $\dim_{K_{t(\rho)}}(M_\rho)=d^\rho_{h(\rho)t(\rho)}$ and $\dim_{K_{h(\rho)}}(M_\rho)=d^\rho_{t(\rho)h(\rho)}$.
  \end{enumerate}

  A \emph{species}\index{Species} (also called a \emph{modulated quiver}\index{Quiver!modulated}) $\Q$ is a pair $(\Delta,\mathbb{M})$, where $\Delta$ is a relative valued quiver and $\mathbb{M}$ is a modulation of $\Delta$.

  A \emph{species morphism} $\Q \rightarrow \Q'$ consists of a relative valued quiver morphism $\varphi: \Delta \rightarrow \Delta'$, a division ring morphism $\psi_i: K_i \rightarrow K'_{\varphi_0(i)}$ for each $i \in \Delta_0$, and a compatible abelian group homomorphism $\psi_\rho: M_\rho \rightarrow M'_{\varphi_1(\rho)}$ for each $\rho \in \Delta_1$. That is, for every $\rho \in \Delta_1$ we have $\psi_\rho(a \cdot m)=\psi_{h(\rho)}(a) \cdot \psi_\rho(m)$ and $\psi_\rho(m \cdot b)=\psi_\rho(m) \cdot \psi_{t(\rho)}(b)$ for all $a \in K_{h(\rho)}$, $b \in K_{t(\rho)}$ and $m \in M_\rho$.

  Let $\spec$\index{$\spec$} denote the category of species.
\end{defn}

\begin{rmk}\label{rmk:multi-arrows} Notice that we allow parallel arrows in our definition of valued quivers and thus in our definition of species. However, many texts only allow for single arrows in their definition of species. We will see in Sections 5 and 6 that we can always assume, without loss of generality, that we have no parallel arrows. Thus our definition of species is consistent with the other definitions in the literature.
\end{rmk}

Another definition of species also appears in the literature (see for example \cite{ddpw}). This definition depends on a central field $K$, and so to distinguish between the two definitions, we will call these objects $K$-species.

\begin{defn}[K-species] 
  Let $\Gamma$ be an absolute valued quiver\index{Quiver!valued!absolute} with valuation \linebreak $(d_i,m_\rho)_{i\in \Gamma_0, \rho \in \Gamma_1}$. A \emph{$K$-modulation}\index{Modulation!$K$-modulation} $\mathbb{M}$ of $\Gamma$ consists of a $K$-division algebra $K_i$ for each $i\in \Gamma_0$, and a $(K_{h(\rho)},K_{t(\rho)})$-bimodule $M_\rho$ for each $\rho \in \Gamma_1$, such that the following two conditions hold:

    \begin{enumerate}
      \item $K$ acts centrally on $M_\rho$ (i.e.\ $k\cdot m=m\cdot k \quad \forall \; k\in K,\; m\in M_\rho$), and
      \item $\dim_K(K_i) = d_i$ and $\dim_K(M_\rho)=m_\rho$.
    \end{enumerate}

  A \emph{$K$-species}\index{Species!$K$-species} (also called a $K$\emph{-modulated quiver}\index{Quiver!$K$-modulated}) $\Q$ is a pair $(\Gamma,\mathbb{M})$, where $\Gamma$ is an absolute valued quiver and $\mathbb{M}$ is a $K$-modulation of $\Gamma$.

  A \emph{$K$-species morphism} $\Q \rightarrow \Q'$ consists of an absolute valued quiver morphism $\varphi: \Gamma \rightarrow \Gamma'$, a $K$-division algebra morphism $\psi_i: K_i \rightarrow K'_{\varphi_0(i)}$ for each $i \in \Gamma_0$, and a compatible $K$-linear map $\psi_\rho: M_\rho \rightarrow M'_{\varphi_1(\rho)}$ for each $\rho \in \Gamma_1$. That is, for every $\rho \in \Gamma_1$ we have $\psi_\rho(a \cdot m)=\psi_{h(\rho)}(a) \cdot \psi_\rho(m)$ and $\psi_\rho(m \cdot b)=\psi_\rho(m) \cdot \psi_{t(\rho)}(b)$ for all $a \in K_{h(\rho)}$, $b \in K_{t(\rho)}$ and $m \in M_\rho$.

  Let $\spec_K$\index{$\spec_K$} denote the category of $K$-species.

\end{defn}

Note that, given a base field $K$, not every absolute valued quiver has a $K$-modulation. For example, it is well-known that the only division algebras over $\mathbb{R}$ are $\mathbb{R}$, $\mathbb{C}$ and $\mathbb{H}$, which have dimension 1, 2 and 4, respectively. Thus, any absolute valued quiver containing a vertex with value 3 (or any value not equal to 1, 2 or 4) has no $\mathbb{R}$-modulation. However, given an absolute valued quiver, we can always find a base field $K$ for which there exists a $K$-modulation. For example, $\mathbb{Q}$ admits field extensions (thus division algebras) of arbitrary dimension, thus $\mathbb{Q}$-modulations always exist.

It is also worth noting that given a valued quiver, relative or absolute, there may exist several non-isomorphic species or $K$-species (depending on the field $K$).

\begin{egg}\label{egg:kspecies}
  Consider the following absolute valued quiver $\Gamma$ and its image under $\F$.

  \begin{center}
    \begin{tikzpicture} [>=stealth]
      \draw (-1,0) node {$\Gamma$:}; \draw (4,0) node {$\F(\Gamma)$:};

      \filldraw (0,0) circle (2pt) (2,0) circle (2pt); \filldraw (5,0) circle (2pt) (7,0) circle (2pt);

      \draw [->] (0,0) -- (1.1,0); \draw (1.1,0) -- (2,0);
      \draw [->] (5,0) -- (6.1,0); \draw (6.1,0) -- (7,0);

      \begin{small}
	\draw (0,-0.3) node {(2)} (2,-0.3) node {(1)} (1,0.3) node {(2)};
	\draw (6,0.3) node {(2,1)};
      \end{small}
    \end{tikzpicture}
  \end{center}
  One can construct the following two $\mathbb{Q}$-species of $\Gamma$.

  \begin{center}
    \begin{tikzpicture} [>=stealth]
      \draw (-1,0) node {$\Q$:} (4,0) node {$\Q'$:};

      \filldraw (0,0) circle (2pt) (2,0) circle (2pt);

      \filldraw (5,0) circle (2pt) (7,0) circle (2pt);

      \draw [->] (0,0) -- (1.1,0); \draw (1.1,0) -- (2,0);

      \draw [->] (5,0) -- (6.1,0); \draw (6.1,0) -- (7,0);

      \begin{small}
	\draw (0,-0.3) node {$\mathbb{Q}(\sqrt{2})$} (1,0.3) node {$\mathbb{Q}(\sqrt{2})$} (2,-0.3) node {$\mathbb{Q}$};

	\draw (5,-0.3) node {$\mathbb{Q}(\sqrt{3})$} (6,0.3) node {$\mathbb{Q}(\sqrt{3})$} (7,-0.3) node {$\mathbb{Q}$};
      \end{small}
    \end{tikzpicture}
  \end{center}
  Then $\Q \ncong \Q'$ as $\mathbb{Q}$-species, since $\mathbb{Q}(\sqrt{2}) \ncong \mathbb{Q}(\sqrt{3})$ as algebras. Also, one can show that $\Q$ and $\Q'$ are species of $\F(\Gamma)$ (indeed, this will follow from Proposition \ref{prop:species}). But again, $\Q \ncong \Q'$ as species since $\mathbb{Q}(\sqrt{2}) \ncong \mathbb{Q}(\sqrt{3})$ as rings.
  
\end{egg}

It is natural to ask how species and $K$-species are related, i.e.\ how the categories $\spec$ and $\spec_K$ are related. To answer this question, we first need the following lemma.

\begin{lem}\label{lem:trace}
  Let $F$ and $G$ be finite-dimensional (nonzero) division algebras over a field $K$ and let $M$ be a finite-dimensional $(F,G)$-bimodule on which $K$ acts centrally. Then $\Hom_F(M,F) \cong \Hom_G(M,G)$ as $(G,F)$-bimodules.
\end{lem}

\begin{proof}
  A proof can be found in \cite[Lemma 3.7]{bn}, albeit using slightly different terminology. For convenience, we present a brief sketch of the proof.

  Let $\tau: F \rightarrow K$ be a nonzero $K$-linear map such that $\tau(ab)=\tau(ba)$ for all $a,b \in F$. Such a map is known to exist; one can take the reduced trace map $F \rightarrow Z(F)$, where $Z(F)$ is the centre of $F$ (see \cite[Chapter IX, Section 2, Proposition 6]{weil}) and compose it with any nonzero map $Z(F) \rightarrow K$. Then $T: \Hom_F(M,F) \rightarrow \Hom_K(M,K)$ defined by $\varphi \mapsto \tau \circ \varphi$ is a $(G,F)$-bimodule isomorphism. By an analogous argument, $\Hom_G(M,G) \cong \Hom_K(M,K)$ completing the proof.
\end{proof}

Given Lemma \ref{lem:trace}, we see that if $\Q$ is a $K$-species with absolute valued quiver $\Gamma$, then $\Q$ is a species with underlying relative valued quiver $\F(\Gamma)$. Also, a $K$-species morphism $\Q \rightarrow \Q'$ is a species morphism when viewing $\Q$ and $\Q'$ as species (because an algebra morphism is a ring morphism and a linear map is a group homomorphism). Thus, we may define a forgetful functor $\U_K: \spec_K \rightarrow \spec$\index{$\U_K$}, which forgets the underlying field $K$ and views absolute valued quivers as relative valued quivers via the functor $\F$\index{$\F$}. This yields the following result.

\begin{prop}\label{prop:species}
  The functor $\U_K: \spec_K \rightarrow \spec$ is faithful and injective on objects. Hence, we may view $\spec_K$ as a subcategory of $\spec$.
\end{prop}

\begin{proof}
  Faithfulness is clear, since $\F$ is faithful and $\U_K$ then simply forgets the underlying field $K$.

  To see that $\U_K$ is injective on objects, suppose $\Q$ and $\Q'$ are $K$-species with modulations $(K_i,M_\rho)_{i \in \Gamma_0, \rho \in \Gamma_1}$ and $(K'_i,M'_\rho)_{i \in \Gamma'_0, \rho \in \Gamma'_1}$, respectively, such that $\U_K(\Q) = \U_K(\Q')$. Then, the underlying (non-valued) quivers of $\Q$ and $\Q'$ are equal. Moreover, $K_i = K'_i$ for all $i \in \Gamma_0=\Gamma'_1$ and $M_\rho = M'_\rho$ for all $\rho \in \Gamma_1 = \Gamma'_1$. So, $\Q=\Q'$ and thus $\U_K$ is injective.
\end{proof}

Note that $\U_K$ is not full (and hence we cannot view $\spec_K$ as a full subcategory of $\spec$) nor is it essentially surjective. In fact, there are objects in $\spec$ which are not of the form $\U_K(\Q)$ for $\Q \in \spec_K$ for any field $K$. The following examples illustrate these points.

\begin{egg}\label{egg:species1}
  Consider $\mathbb{C}$ as a $\mathbb{C}$-species, that is $\mathbb{C}$ is a $\mathbb{C}$-modulation of the trivially valued quiver consisting of one vertex and no arrows. Then, the only morphism in $\Hom_\mathbb{C}(\mathbb{C},\mathbb{C})$ is the identity, since any such morphism must send $1$ to $1$ and be $\mathbb{C}$-linear. However, $\Hom_{\U_\mathbb{C}(\mathbb{C})} (\U_\mathbb{C}(\mathbb{C}), \U_\mathbb{C}(\mathbb{C}))$ contains more than just the identity. Indeed, let $\varphi:\mathbb{C} \rightarrow \mathbb{C}$ given by $z \mapsto \overline{z}$. Then, $\varphi$ is a ring morphism and thus defines a species morphism. Hence,
  \[ \U_\mathbb{C} : \Hom_\mathbb{C}(\mathbb{C},\mathbb{C}) \rightarrow \Hom_{\U_\mathbb{C}(\mathbb{C})} (\U_\mathbb{C}(\mathbb{C}), \U_\mathbb{C}(\mathbb{C})) \]
  is not surjective and so $\U_\mathbb{C}$ is not full.
\end{egg}

\begin{egg}\label{egg:species2}
  There exist division rings which are not finite-dimensional over their centres; such division rings are called \emph{centrally infinite}. Hilbert was the first to construct such a ring (see for example \cite[Proposition 14.2]{tyl}). Suppose $R$ is a centrally infinite ring. Then, for any field $K$ contained in $R$ such that $R$ is a $K$-algebra, $K \subseteq Z(R)$ and so $R$ is not finite-dimensional $K$-algebra. Thus, any species containing $R$ as part of its modulation is not isomorphic to any object in the image of $\U_K$ for any field $K$.
\end{egg}

One might think that we could eliminate this problem by restricting ourselves to modulations containing only centrally finite rings. In other words, one might believe that if $\Q$ is a species whose modulation contains only centrally finite rings, then we can find a field $K$ and a $K$-species $\Q'$ such that $\Q \cong U_K(\Q')$. However, this is not the case as we see in the following example.

\begin{egg}\label{egg:species3}
  Let $p$ be a prime. Consider:

  \begin{center}
      \begin{tikzpicture} [>=stealth, baseline=0]

      \draw (-1,0) node {$\Q$:};

      \filldraw (0,0) circle (2pt) (1,0) circle (2pt);

      \draw [->] (0,0) -- (0.6,0); \draw (0.6,0) -- (1,0);

      \begin{small}
	\draw (0,-0.3) node {$G$} (1,-0.3) node {$F$} (0.5,0.3) node {$M$};
      \end{small}
    \end{tikzpicture}$\;$,
  \end{center}
  where $F=G=\overline{\mathbb{F}}_p$ ($F$ and $G$ are then centrally finite since they are fields) and $M=\overline{\mathbb{F}}_p$ is an $(F,G)$-bimodule with actions:
  \[ f \cdot m \cdot g = fmg^p\]
  for all $f \in F$, $g \in G$ and $m \in M$. We claim that $\Q$ is a species. The dimension criterion is clear, as $\dim_F M = \dim_G M = 1$. Thus, it remains to show that
  \[ \Hom_F(M,F) \cong \Hom_G(M,G)\]
  as $(G,F)$-bimodules. Recall that in $\overline{\mathbb{F}}_p$, $p$-th roots exist and are unique. Indeed, for any $a \in \overline{\mathbb{F}}_p$, the $p$-th roots of $a$ are the roots of the polynomial $x^p - a$. Because $\overline{\mathbb{F}}_p$ is algebraically closed, this polynomial has a root, say $\alpha$. Because $\ch \overline{\mathbb{F}}_p = p$ we have
  \[ (x-\alpha)^p = x^p - \alpha^p = x^p - a. \]
  Hence, $\alpha$ is the unique $p$-th root of $a$. Therefore, we have a well-defined map:
  \begin{align*}
    \Phi: \Hom_F(M,F) & \rightarrow \Hom_G(M,G) \\
    \varphi & \mapsto \rho \circ \varphi,
  \end{align*}
  where $\rho$ is the $p$-th root map. It is straightforward to show that $\Phi$ is a $(G,F)$-bimodule isomorphism.

  Therefore, $\Q$ is a species. Yet the field $\overline{\mathbb{F}}_p$ does not act centrally on $M$. Indeed, take an element $a \notin \mathbb{F}_p$, then
  \[ a \cdot 1 = a \ne a^p = 1 \cdot a. \]
  In fact, the only subfield that does act centrally on $M$ is $\mathbb{F}_p$ since $a^p = a$ if and only if $a \in \mathbb{F}_p$. But, $F$ and $G$ are infinite-dimensional over $\mathbb{F}_p$. Thus, there is no field $K$ for which $\Q$ is isomorphic to an object of the form $\U_K(\Q')$ with $\Q' \in \spec_K$.
\end{egg}

\section{The Path and Tensor Algebras}
In this section we will define the path and tensor algebras associated to quivers and species, respectively. These algebras play an important role in the representation theory of finite-dimensional algebras (see Theorems \ref{theo:pathalg} and \ref{theo:tensalg}, and Corollary \ref{cor:alg}). In subsequent sections, we will give a more in-depth study of these algebras (Sections 4 and 5) and we will show that modules of path and tensor algebras are equivalent to representations of quivers and species, respectively (Section 6).

Recall that a \emph{path}\index{Path} of length $n$ in a quiver\index{Quiver} $Q$ is a sequence of $n$ arrows in $Q_1$, $\rho_n \rho_{n-1} \cdots \rho_1$, such that $h(\rho_i)=t(\rho_{i+1})$ for all $i=1,\;2, \dots, \; n-1$. For every vertex, we have a trivial path of length 0 (beginning and ending at that vertex).

\begin{defn}[Path algebra] 
  The \emph{path algebra}\index{Path algebra}, $KQ$\index{$KQ$}, of a quiver $Q$ is the $K$-algebra with basis the set of all the paths in $Q$ and multiplication given by:

  \begin{center}
    $(\beta_n \beta_{n-1} \cdots \beta_1)(\alpha_m \alpha_{m-1} \cdots \alpha_1) = 
    \begin{cases}
      \beta_n \beta_{n-1} \cdots \beta_1\alpha_m \alpha_{m-1} \cdots \alpha_1, & \text{if } t(\beta_1)=h(\alpha_m), \\
      0,								       & \text{otherwise}.
    \end{cases}$

  \end{center}
\end{defn}

\begin{rmk}
  According to the convention used, a path $i_1 \xrightarrow{\rho_1} i_2 \xrightarrow{\rho_2} \cdots \xrightarrow{\rho_{n-1}} i_n \xrightarrow{\rho_n} i_{n+1}$ is written from ``right to left'' $\rho_n \rho_{n-1} \cdots \rho_1$. However, some texts write paths from ``left to right'' $\rho_1 \rho_2 \cdots \rho_n$. Using the ``left to right'' convention yields a path algebra that is opposite to the one defined here.
\end{rmk}

Note that $KQ$ is associative and unital (its identity is $\sum_{i \in Q_0} \varepsilon_i$, where $\varepsilon_i$ is the path of length zero at $i$). Also, $KQ$ is finite-dimensional precisely when $Q$ contains no oriented cycles.

\begin{defn}[Admissible ideal] 
  Let $Q$ be a quiver and let $P^n(Q)=\vspan_K$\{all paths in $Q$ of length $ \ge n$\}. An \emph{admissible ideal}\index{Admissible ideal} $I$ of the path algebra $KQ$ is a two-sided ideal of $KQ$ satisfying

  \begin{center}
    $P^n(Q) \subseteq I \subseteq P^2(Q), \quad$ for some positive integer $n$.
  \end{center}
\end{defn}

If $Q$ has no oriented cycles, then any ideal $I \subseteq P^2(Q)$ of $KQ$ is an admissible ideal, since $P^n(Q)=0$ for sufficiently large $n$.

There is a strong relationship between path algebras and finite-dimensional algebras, touched upon by Brauer \cite{brauer}, Jans \cite{jans} and Yoshii \cite{yoshii}, but fully explored by Gabriel \cite{g}. Let $A$ be a finite-dimensional $K$-algebra. We recall a few definitions. An element $\varepsilon \in A$ is called an \emph{idempotent}\index{Idempotent} if $\varepsilon^2 = \varepsilon$. Two idempotents $\varepsilon_1$ and $\varepsilon_2$ are called \emph{orthogonal}\index{Idempotent!orthogonal} if $\varepsilon_1 \varepsilon_2 = \varepsilon_2 \varepsilon_1 = 0$. An idempotent $\varepsilon$ is called \emph{primitive}\index{Idempotent!primitive} if it cannot be written as a sum $\varepsilon = \varepsilon_1 + \varepsilon_2$, where $\varepsilon_1$ and $\varepsilon_2$ are orthogonal idempotents. A set of idempotents $\{\varepsilon_1, \dots, \varepsilon_n\}$ is called \emph{complete}\index{Idempotent!complete set of} if $\sum_{i=1}^n \varepsilon_i=1$. If $\{\varepsilon_1, \dots, \varepsilon_n\}$ is a complete set of primitive (pairwise) orthogonal idempotents of $A$, then $A=A\varepsilon_1 \oplus \dots \oplus A\varepsilon_n$ is a decomposition of $A$ (as a left $A$-module) into indecomposable modules; this decomposition is unique up to isomorphism and permutation of the terms. We say that $A$ is \emph{basic}\index{Basic algebra} if $A\varepsilon_i \ncong A\varepsilon_j$ as (left) $A$-modules for all $i \ne j$ (or, alternatively, the decomposition of $A$ into indecomposable modules admits no repeated factors). Finally, $A$ is called \emph{hereditary}\index{Hereditary algebra} if every $A$-submodule of a projective $A$-module is again projective.

\begin{theo}\label{theo:pathalg} Let $K$ be an algebraically closed field and let $A$ be a finite-dimensional $K$-algebra.

  \begin{enumerate}
    \item If $A$ is basic and hereditary, then $A \cong KQ$ (as $K$-algebras) for some quiver $Q$.
    \item If $A$ is basic, then $A \cong KQ/I$ (as $K$-algebras) for some quiver $Q$ and some admissible ideal $I$ of $KQ$.
  \end{enumerate}

\end{theo}

For a proof of Theorem \ref{theo:pathalg} see \cite[Sections II and VII]{ass} or \cite[Propositions 4.1.7 and 4.2.4]{b} (though it also follows from \cite[Proposition 10.2]{dr2}). The above result is powerful, but it does not necessarily hold over fields which are not algebraically closed. If we want to work with algebras over non-algebraically closed fields, we need to generalize the notion of a path algebra. We look, then, at the analogue of the path algebra for a $K$-species\index{Species!$K$-species}.

Let $\Q$ be a species\index{Species} of a relative valued quiver $\Delta$ with modulation $(K_i,M_\rho)_{i \in \Delta_0, \rho \in \Delta_1}$. Let $D=\Pi_{i\in \Delta_0}K_i$ and let $M=\bigoplus_{\rho \in \Delta_1} M_\rho$. Then $D$ is a ring and $M$ naturally becomes a $(D,D)$-bimodule. If $\Q$ is a $K$-species, then $D$ is a $K$-algebra.

\begin{defn}[Tensor ring/algebra] 
  The \emph{tensor ring}\index{Tensor ring}, $T(\Q)$\index{$T(\Q)$}, of a species $\Q$ is defined by
  \[ T(\mathcal{Q})=\bigoplus_{n=0}^\infty T^n(M), \]
  where
  \[ T^0(M)=D \; \text{ and } \; T^n(M)=T^{n-1}(M)\otimes_D M  \text{ for } n \ge 1. \]
  Multiplication is determined by the composition
  \[ T^m(M) \times T^n(M) \twoheadrightarrow T^m(M) \otimes_D T^n(M) \xrightarrow{\cong} T^{m+n}(M). \]
  If $\Q$ is a $K$-species, then $T(\Q)$ is a $K$-algebra. In this case we call $T(\Q)$ the \emph{tensor algebra}\index{Tensor algebra} of $\Q$.
\end{defn}

Admissible ideals\index{Admissible ideal} for tensor rings/algebras are defined in the same way as admissible ideals for path algebras by setting $P^n(\Q)=\bigoplus_{m=n}^\infty T^m(M)$.

Suppose that $\Gamma$ is an absolute valued quiver with trivial valuation (all $d_i$ and $m_\rho$ are equal to 1) and $\Q$ is a $K$-species of $\Gamma$. Then, for each $i \in \Gamma_0$, $\dim_K K_i=1$, which implies that $K_i \cong K$ (as $K$-algebras). Likewise, $\dim_K M_\rho = 1$ implies that $M_\rho \cong K$ (as $(K,K)$-bimodules). Therefore, it follows that $T(\Q) \cong KQ$ where $Q=(\Gamma_0,\Gamma_1)$. Thus, when viewing non-valued quivers as absolute valued quivers with trivial valuation, the tensor algebra of the $K$-species becomes simply the path algebra (over $K$) of the quiver. Therefore, the tensor algebra is indeed a generalization of the path algebra. Additionally, the tensor algebra allows us to generalize Theorem \ref{theo:pathalg}.

Recall that a field $K$ is called \emph{perfect}\index{Perfect field} if either $\ch(K)=0$ or, if $\ch(K)=p > 0$, then $K^p=\{a^p \; | \; a \in K\}=K$.

\begin{theo}\label{theo:tensalg} Let $K$ be a perfect field and let $A$ be a finite-dimensional $K$-algebra.
 
  \begin{enumerate}
   \item If $A$ is basic and hereditary, then $A \cong T(\Q)$ (as $K$-algebras) for some $K$-species $\Q$.
   \item If $A$ is basic, then $A \cong T(\Q)/I$ (as $K$-algebras) for some $K$-species $\Q$ and some admissible ideal $I$ of $T(\Q)$.
  \end{enumerate}

\end{theo}

For a proof of Theorem \ref{theo:tensalg}, see \cite[Proposition 10.2]{dr2} or \cite[Corollary 4.1.11 and Proposition 4.2.5]{b} or \cite[Section 8.5]{dk}. Note that Theorem \ref{theo:tensalg} does not necessarily hold over non-perfect fields. To see why, we first introduce a useful tool in the study of path and tensor algebras.

\begin{defn}[Jacobson radical] 
  The \emph{Jacobson radical}\index{Jacobson radical} of a ring $R$ is the intersection of all maximal left ideals of $R$. We denote the Jacobson radical of $R$ by $\rad R$ \index{$\rad$}.
\end{defn}

\begin{rmk}
  The intersection of all maximal left ideals coincides with the intersection of all maximal right ideals (see, for example, \cite[Corollary 4.5]{tyl}), so the Jacobson radical could alternatively be defined in terms of right ideals.
\end{rmk}

\begin{lem}\label{lem:rad}
  Let $\Q$ be a species.

  \begin{enumerate}
    \item If $\Q$ contains no oriented cycles, then $\rad T(\Q) = \bigoplus_{n=1}^\infty T^n(M)$.
    \item Let $I$ be an admissible ideal of $T(\Q)$. Then, $\rad \left( T(\Q)/I \right) = \left( \bigoplus_{n=1}^\infty T^n(M) \right) /I$.
  \end{enumerate}

\end{lem}

\begin{proof}
  It is well known that if $R$ is a ring and $J$ is a two-sided nilpotent ideal of $R$ such that $R/J$ is semisimple, then $\rad R = J$ (see for example \cite[Lemma 4.11 and Proposition 4.6]{tyl} together with the fact that the radical of a semisimple ring is 0). Let $J= \bigoplus_{n=1}^\infty T^n(M)$. If $\Q$ contains no oriented cycles, then $T^n(M)=0$ for some positive integer $n$. Thus, $J^n=0$ and $J$ is then nilpotent. Then $T(\Q)/J \cong T^0(M) = D$, which is semisimple. Therefore, $\rad T(\Q) = J$, proving Part 1.

  If $I$ is an admissible ideal of $T(\Q)$, let $J=\left( \bigoplus_{n=1}^\infty T^n(M)/I \right)$. By definition, $P^n(\Q) \subseteq I$ for some $n$ and so $J^n=0$. Thus, Part 2 follows by a similar argument.
\end{proof}

\begin{rmk}
  Part 1 of Lemma \ref{lem:rad} is false if $\Q$ contains oriented cycles. One does not need to look beyond quivers to see why. For example, following \cite[Section II, Chapter 1]{ass}, we can consider the path algebra of the Jordan quiver over an infinite field $K$. That is, we consider $KQ$, where:

  
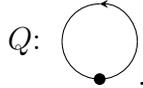
\begin{center}
    \begin{tikzpicture}[>=stealth]
      \draw (-1,0.5) node {$Q$:};

      \filldraw (0,0) circle (2pt);

      \draw [out=0, in=-90] (0,0) to (0.5,0.5); \draw[->] [out=90, in=0] (0.5,0.5) to (0,1);
      \draw [out=180, in=90] (0,1) to (-0.5,0.5); \draw [out=-90, in=180] (-0.5,0.5) to (0,0);
    \end{tikzpicture}.
  \end{center}
  Then it is clear that $KQ \cong K[t]$, the polynomial ring in one variable. For each $\alpha \in K$, let $I_\alpha$ be the ideal generated by $t+\alpha$. Each $I_\alpha$ is a maximal ideal and $\bigcap_{\alpha \in K} I_\alpha = 0$ since $K$ is infinite. Thus $\rad KQ = 0$ whereas $\bigoplus_{n=1}^\infty T^n(Q) \cong (t)$ (the ideal generated by the lone arrow of $Q$).
\end{rmk}

With the concept of the Jacobson radical and Lemma \ref{lem:rad}, we are ready to see why Theorem \ref{theo:tensalg} fails over non-perfect fields. Recall that a $K$-algebra epimorphism $\varphi:A \twoheadrightarrow B$ is said to \emph{split}\index{Split epimorphism} if there exists a $K$-algebra morphism $\mu:B \rightarrow A$ such that $\varphi \circ \mu = \id_B$. We see that if $A = T(\Q) / I$ for a $K$-species $\Q$ and admissible ideal $I$, then the canonical projection $A \twoheadrightarrow A/\rad A$ splits (since $A \cong D \oplus \rad A$). Thus, to construct an example where Theorem \ref{theo:tensalg} fails, it suffices to find an algebra where this canonical projection does not split. This is possible over a non-perfect field.

\begin{egg} \cite[Remark (ii) following Corollary 4.1.11]{b}
  Let $K_0$ be a field of characteristic $p \ne 0$ and let $K=K_0(t)$, which is not a perfect field. Let $A=K[x,y]/(x^p,y^p-x-t)$. A quick calculation shows that $\rad A = (x)$ and thus $A / \rad A \cong K[y]/(y^p-t)$. One can easily verify that the projection $A \twoheadrightarrow A/ \rad A$ does not split. Hence, $A$ is not isomorphic to the quotient of the tensor algebra of a species by some admissible ideal.
\end{egg}

Theorems \ref{theo:pathalg} and \ref{theo:tensalg} require our algebras to be basic\index{Basic algebra}. There is a slightly weaker property that holds in the case of non-basic algebras, that of \emph{Morita equivalence}\index{Morita equivalence}.

\begin{defn}[Morita equivalence]
  Two rings $R$ and $S$ are said to be \emph{Morita equivalent}\index{Morita equivalence} if their categories of (left) modules, $R$-Mod and $S$-Mod, are equivalent.
\end{defn}

\begin{cor}\label{cor:alg} Let $K$ be a field and let $A$ be a finite-dimensional $K$-algebra.
  \begin{enumerate}
    \item If $K$ is algebraically closed and $A$ is hereditary, then $A$ is Morita equivalent to $KQ$ for some quiver $Q$.
    \item If $K$ is algebraically closed, then $A$ is Morita equivalent to $KQ/I$ for some quiver $Q$ and some admissible ideal $I$ of $KQ$.
    \item If $K$ is perfect and $A$ is hereditary, then $A$ is Morita equivalent to $T(\Q)$ for some $K$-species $\Q$.
    \item If $K$ is perfect, then $A$ is Morita equivalent to $T(\Q)/I$ for some $K$-species $\Q$ and some admissible ideal $I$ of $T(\Q)$.
  \end{enumerate}
\end{cor}

\begin{proof}
  Every algebra is Morita equivalent to a basic algebra (see \cite[Section 2.2]{b}) and Morita equivalence preserves the property of being hereditary (indeed, an equivalence of categories preserves projective modules). Thus, the result follows as a consequence of Theorems \ref{theo:pathalg} and \ref{theo:tensalg}.
\end{proof}

\section{The Frobenius Morphism}
When working over the finite field of $q$ elements, $\mathbb{F}_q$, it is possible to avoid dealing with species altogether and deal only with quivers with automorphism. This is achieved by using the Frobenius morphism (described below).

\begin{defn}[Frobenius morphism] 
  Let $K=\overline{\mathbb{F}}_q$ (the algebraic closure of $\mathbb{F}_q$). Given a quiver with automorphism $(Q,\sigma)$, the \emph{Frobenius morphism}\index{Frobenius morphism} $F=F_{Q,\sigma,q}$\index{$F$} is defined as
  \begin{align*}
    F:KQ 	 &\rightarrow KQ \\
    \sum_i \lambda_i p_i &\mapsto \sum_i \lambda_i^q \sigma(p_i)
  \end{align*}
  for all $\lambda_i \in K$ and paths $p_i$ in $Q$. The $F$-fixed point algebra is
  \[ (KQ)^F=\{x \in KQ \;|\; F(x)=x\}. \]
\end{defn}

Note that while $KQ$ is an algebra over $K$, the fixed point algebra $(KQ)^F$ is an algebra over $\mathbb{F}_q$. Indeed, suppose $0 \ne x \in (KQ)^F$, then $F(\lambda x)=\lambda^q F(x)=\lambda^q x$. Thus, $\lambda x \in (KQ)^F$ if and only if $\lambda^q=\lambda$, which occurs if and only if $\lambda \in \mathbb{F}_q$.

Suppose $\Gamma$ is the absolute valued quiver obtained by folding\index{Folding} $(Q,\sigma)$. For each $i \in \Gamma_0$ and each $\rho \in \Gamma_1$ define
\[ A_i=\bigoplus_{a \in i} K\varepsilon_a \text{ and } A_\rho=\bigoplus_{\tau \in \rho} K\tau. \]
where $\varepsilon_a$ is the trivial path at vertex $a$. Then as an $\mathbb{F}_q$-algebra, $A_i^F \cong \mathbb{F}_{q^{d_i}}$. Indeed, fix some $a \in i$, then,
\[ A_i^F = \left\{ \left. x=\sum_{j=0}^{d_i-1}\lambda_j \varepsilon_{\sigma^j(a)} \; \right| \; \lambda_j \in K \text{ and } F(x)=x \right\} \]
Applying $F$ to an arbitrary $x=\sum_{j=0}^{d_i-1}\lambda_j \varepsilon_{\sigma^j(a)} \in A_i^F$, we obtain:
\[ F(x) = \sum_{j=0}^{d_i-1}\lambda_j^q \sigma(\varepsilon_{\sigma^j(a)})=\sum_{j=0}^{d_i-1}\lambda_j^q \varepsilon_{\sigma^{j+1}(a)}. \]
The equality $F(x)=x$ yields $\lambda_j^q=\lambda_{j+1}$ for $j=0,1, \dots, d_i-2$ and $\lambda_{d_i-1}^q=\lambda_0$. By successive substitution, we get $\lambda_0^{q^{d_i}}=\lambda_0$, which occurs if and only if $\lambda_0 \in \mathbb{F}_{q^{d_i}}$, and $\lambda_j=\lambda_0^{q^j}$. Thus, $A_i^F$ can be rewritten as:
\begin{align*}
  A_i^F &= \left\{ \left. \sum_{j=0}^{d_i-1} \lambda_0^{q^j} \varepsilon_{\sigma^j(a)} \; \right| \; \lambda_0 \in \mathbb{F}_{q^{d_i}} \right\} \\
	&\cong \mathbb{F}_{q^{d_i}} \qquad \text{(as fields).}
\end{align*}

It is easy to see that $A_\rho$ is an $(A_{h(\rho)},A_{t(\rho)})$-bimodule via multiplication, thus $A_\rho^F$ is an $(A_{h(\rho)}^F,A_{t(\rho)}^F)$-bimodule (on which $\mathbb{F}_q$ acts centrally). Since $A_i^F \cong \mathbb{F}_{q^{d_i}}$ for each $i \in \Gamma_0$, $A_\rho^F$ is then an $(\mathbb{F}_{q^{d_{h(\rho)}}},\mathbb{F}_{q^{d_{t(\rho)}}})$-bimodule. Over fields, we make no distinction between left and right modules because of commutativity. Thus, $A_\rho^F$ is an $\mathbb{F}_{q^{d_{h(\rho)}}}$-module and an $\mathbb{F}_{q^{d_{t(\rho)}}}$-module, and hence $A_\rho^F$ is a module of the composite field of $\mathbb{F}_{q^{d_{h(\rho)}}}$ and $\mathbb{F}_{q^{d_{t(\rho)}}}$, which in this case is simply the bigger of the two fields (recall that the composite of two fields is the smallest field containing both fields). Over fields, all modules are free and thus $A_\rho^F$ is a free module of the composite field (this fact will be useful later on). Also, $\dim_{\mathbb{F}_q}A_i^F=d_i$ and $\dim_{\mathbb{F}_q}A_\rho^F=m_\rho$ (the dimensions are the number of vertices/arrows in the corresponding orbits). Therefore, $\mathbb{M}=(A_i^F,A_\rho^F)_{i \in \Gamma_0, \rho \in \Gamma_1}$ defines an $\mathbb{F}_q$-modulation of $\Gamma$. We will denote the $\mathbb{F}_q$-species ($\Gamma$,$\mathbb{M}$) by $\mathcal{Q}_{Q,\sigma,q}$. This leads to the following result.

\begin{theo}\textbf{\emph{\cite[Theorem 3.25]{ddpw}}}\label{theo:frob}
  Let $(Q,\sigma)$ be a quiver with automorphism. Then $(KQ)^F \cong T(\mathcal{Q}_{Q,\sigma,q})$ as $\mathbb{F}_q$-algebras\index{$KQ$}\index{Path algebra}\index{$T(\Q)$}\index{Tensor algebra}.
\end{theo}

In light of Theorem \ref{theo:frob}, the natural question to ask is: given an arbitrary $\mathbb{F}_q$-species, is its tensor algebra isomorphic to the fixed point algebra of a quiver with automorphism? And if so, to which one?

Suppose $\Q$ is an $\mathbb{F}_q$-species with underlying absolute valued quiver $\Gamma$ and $\mathbb{F}_q$-modulation $(K_i,M_\rho)_{i \in \Gamma_0, \rho \in \Gamma_1}$. Each $K_i$ is, by definition, a division algebra containing $q^{d_i}$ elements. According to the well-known Wedderburn's little theorem, all finite division algebras are fields. Thus, $K_i \cong \mathbb{F}_{q^{d_i}}$. Similar to the above discussion, $M_\rho$ is then a free module of the composite field of $\mathbb{F}_{q^{d_{h(\rho)}}}$ and $\mathbb{F}_{q^{d_{t(\rho)}}}$. Therefore, by unfolding $\Gamma$ (as in Section 1, say) to a quiver with automorphism $(Q,\sigma)$, we get $\Q \cong \Q_{Q,\sigma,q}$ as $\mathbb{F}_q$-species. This leads to the following result.

\begin{prop}\textbf{\emph{\cite[Proposition 3.37]{ddpw}}}
 For any $\mathbb{F}_q$-species $\mathcal{Q}$, there exists a quiver with automorphism $(Q,\sigma)$ such that $T(\mathcal{Q}) \cong (KQ)^F$ as $\mathbb{F}_q$-algebras.
\end{prop}

Note that, given an $\mathbb{F}_q$-species $\Q=(\Gamma,\mathbb{M})$ and a quiver with automorphism $(Q,\sigma)$ such that $T(\mathcal{Q}) \cong (KQ)^F$, we cannot to conclude that $(Q,\sigma)$ folds into $\Gamma$ as the following example illustrates.

\begin{egg} Consider the following quiver.

  \begin{center}
    \begin{tikzpicture} [>=stealth]
      \draw (-1,0) node {$Q$:};

      \filldraw (0,0) circle (2pt) (1,0) circle (2pt);

      \draw [->] (0,0.05) -- (0.6,0.05); \draw (0.5,0.05) -- (1,0.05);
      \draw [->] (0,-0.05) -- (0.6,-0.05); \draw (0.5,-0.05) -- (1,-0.05);
    \end{tikzpicture}
  \end{center}

  There are two possible automorphisms of $Q$: $\sigma=\id_Q$ and $\sigma'$, the automorphism defined by interchanging the two arrows of $Q$. By folding $Q$ with respect to $\sigma$ and $\sigma'$, we obtain the following two absolute valued quivers.

  \begin{center}
    \begin{tikzpicture} [>=stealth]
      \draw (-1,0) node {$\Gamma$:} (4,0) node {$\Gamma'$:};

      \filldraw (0,0) circle (2pt) (1,0) circle (2pt);

      \filldraw (5,0) circle (2pt) (6,0) circle (2pt);

      \draw [->] (0,0.05) -- (0.6,0.05); \draw (0.5,0.05) -- (1,0.05);
      \draw [->] (0,-0.05) -- (0.6,-0.05); \draw (0.5,-0.05) -- (1,-0.05);

      \draw [->] (5,0) -- (5.6,0); \draw (5.5,0) -- (6,0);

      \begin{small}
	\draw (-0.3,0) node {(1)} (0.5,0.3) node {(1)} (0.5,-0.3) node {(1)} (1.3,0) node {(1)};

	\draw (4.7,0) node {(1)} (5.5,0.3) node {(2)} (6.3,0) node {(1)};
      \end{small}
    \end{tikzpicture}
  \end{center}

  It is clear that $\Gamma$ and $\Gamma'$ are not isomorphic. Now, construct $\mathbb{F}_q$-species of $\Gamma$ and $\Gamma'$ with the following $\mathbb{F}_q$-modulations.

  \begin{center}
    \begin{tikzpicture} [>=stealth]
      \draw (-1,0) node {$\mathcal{Q}$:} (4,0) node {$\mathcal{Q}'$:};

      \filldraw (0,0) circle (2pt) (1,0) circle (2pt);

      \filldraw (5,0) circle (2pt) (6,0) circle (2pt);

      \draw [->] (0,0.05) -- (0.6,0.05); \draw (0.5,0.05) -- (1,0.05);
      \draw [->] (0,-0.05) -- (0.6,-0.05); \draw (0.5,-0.05) -- (1,-0.05);

      \draw [->] (5,0) -- (5.6,0); \draw (5.5,0) -- (6,0);

      \begin{small}
	\draw (-0.3,0) node {$\mathbb{F}_q$} (0.5,0.3) node {$\mathbb{F}_q$} (0.5,-0.3) node {$\mathbb{F}_q$} (1.3,0) node {$\mathbb{F}_q$};

	\draw (4.7,0) node {$\mathbb{F}_q$} (5.5,0.3) node {$\mathbb{F}_q^2$} (6.3,0) node {$\mathbb{F}_q$};
      \end{small}
    \end{tikzpicture}
  \end{center}

  Then we have that $T(\mathcal{Q}) \cong T(\mathcal{Q}')$. Thus, $(KQ)^{F_{Q,\sigma,q}} \cong T(\Q')$, yet $(Q,\sigma)$ does not fold into $\Gamma'$.

\end{egg}

The above example raises an interesting question. Notice that the two $\mathbb{F}_q$-species $\Q$ and $\Q'$ are not isomorphic, but their tensor algebras $T(\Q)$ and $T(\Q')$ are isomorphic. This phenomenon is not restricted to finite fields either; if we replaced $\mathbb{F}_q$ with some arbitrary field $K$, we still get $\Q \ncong \Q'$ as $K$-species, but $T(\Q) \cong T(\Q')$ as $K$-algebras. So we may ask: under what conditions are the tensor algebras (or rings) of two $K$-species (or species) isomorphic? We answer this question in the following section.

It is worth noting that over infinite fields, there are no (known) methods to extend the results of this section. It is tempting to think that, given an infinite field $K$ and a quiver with automorphism $(Q,\sigma)$ that folds into an absolute valued quiver $\Gamma$, Theorem \ref{theo:frob} might be extended by saying that the fixed point algebra $(KQ)^\sigma$ ($\sigma$ extends to an automorphism of $KQ$) is isomorphic to the tensor algebra of a $K$-species of $\Gamma$. This is, however, not the case as the following example illustrates.

\begin{egg}
  Take $K=\mathbb{R}$ to be our (infinite) base field. Consider the following quiver.

  \begin{center}
    \begin{tikzpicture}[>=stealth]
      \draw (-1,0) node {$Q:$};

      \filldraw (0,0) circle (2pt) (1,0) circle (2pt) (2,0) circle (2pt);

      \draw[->] (0,0) -- (0.6,0); \draw (0.6,0) -- (1,0);
      \draw[->] (2,0) -- (1.4,0); \draw (1.4,0) -- (1,0);

      \draw (0,-0.3) node {$\varepsilon_1$} (1,-0.3) node {$\varepsilon_2$} (2,-0.3) node {$\varepsilon_3$};
      \draw (0.5,0.3) node {$\alpha$} (1.5,0.3) node {$\beta$};
    \end{tikzpicture}
  \end{center}
  Let $\sigma$ be the automorphism of $Q$ given by:
  \[ \sigma :\begin{pmatrix}
      \varepsilon_1 & \varepsilon_2 & \varepsilon_3 & \alpha & \beta \\
      \varepsilon_3 & \varepsilon_2 & \varepsilon_1 & \beta & \alpha
     \end{pmatrix}. \]
  Then $(Q,\sigma)$ folds into the following absolute valued quiver.

  \begin{center}
    \begin{tikzpicture}[>=stealth]
      \draw (-1,0) node {$\Gamma:$};

      \filldraw (0,0) circle (2pt) (1,0) circle (2pt);

      \draw[->] (0,0) -- (0.6,0); \draw (0.6,0) -- (1,0);

      \draw (0,-0.3) node {$(2)$} (0.5,0.3) node {$(2)$} (1,-0.3) node {$(1)$};
    \end{tikzpicture}
  \end{center}
  So, we would like for $(KQ)^\sigma$ to be isomorphic to the tensor algebra of a $K$-species of $\Gamma$. However, this does not happen.

  An element $x=a_1 \varepsilon_1 + a_2 \varepsilon_2 + a_3 \varepsilon_3 + a_4 \alpha + a_5 \beta \in KQ$ is fixed by $\sigma$ if and only if $\sigma(x) = x$, that is, if and only if
  \[ a_3 \varepsilon_1 + a_2 \varepsilon_2 + a_1 \varepsilon_1 + a_5 \alpha + a_4 \beta = a_1 \varepsilon_1 + a_2 \varepsilon_2 + a_3 \varepsilon_3 + a_4 \alpha + a_5 \beta, \]
  which occurs if and only if $a_1 = a_3$ and $a_4 = a_5$. Hence, $(KQ)^\sigma$ has basis $\{ \varepsilon_1 + \varepsilon_3, \varepsilon_2, \alpha+\beta \}$ and we see that it is isomorphic to the path algebra (over $K$) of the following quiver.

  \begin{center}
    \begin{tikzpicture}[>=stealth]
      \draw (-1,0) node {$Q':$};

      \filldraw (0,0) circle (2pt) (1.5,0) circle (2pt);

      \draw[->] (0,0) -- (0.85,0); \draw (0.85,0) -- (1.5,0);

      \draw (0,-0.3) node {$\varepsilon_1 + \varepsilon_3$} (0.75,0.3) node {$\alpha + \beta$} (1.5,-0.3) node {$\varepsilon_2$};
    \end{tikzpicture}
  \end{center}
  The algebra $KQ'$ is certainly not isomorphic to the tensor algebra of any $K$-species of $\Gamma$. Indeed, over $K$, $KQ'$ has dimension $3$ whereas the tensor algebra of any $K$-species of $\Gamma$ has dimension $5$.
\end{egg}

\section{A Closer Look at Tensor Rings}
In this section, we find necessary and sufficient conditions for two tensor \linebreak rings/algebras to be isomorphic. We show that the isomorphism of tensor \linebreak rings/algebras corresponds to an equivalence on the level of species (see Theorem \ref{theo:iso}). Furthermore, we show that if $\Q$ is a $K$-species, where $K$ is a perfect field, then $\overline{K} \otimes_K T(\Q)$ is either isomorphic to, or Morita equivalent to, the path algebra of a quiver (see Theorem \ref{theo:tensor} and Corollary \ref{cor:tensor}). This serves as a partial generalization to \cite[Lemma 21]{hub} in which Hubery proved a similar result when $K$ is a finite field. We begin by introducing the notion of ``crushing''.

\begin{defn}[Crushed absolute valued quiver] 
  Let $\Gamma$ be an absolute valued quiver\index{Quiver!valued!absolute}. Define a new absolute valued quiver, which we will denote $\Gamma^C$, as follows:

  \begin{itemize}
    \item $\Gamma_0^C = \Gamma_0$,
    \item \# arrows $i \rightarrow j = 
            \begin{cases}
              1, & \text{if } \exists \; \rho:i \rightarrow j \in \Gamma_1, \\
              0, & \text{otherwise},
            \end{cases}$
    \item $d_i^C = d_i$ for all $i \in \Gamma_0^C=\Gamma_0$,
    \item $m_\rho^C = \sum_{(\alpha:t(\rho)\rightarrow h(\rho)) \; \in \; \Gamma_1} m_\alpha$ for all $\rho \in \Gamma_1^C$.
  \end{itemize}

  Intuitively, one ``crushes'' all parallel arrows of $\Gamma$ into a single arrow and sums up the values. 

  
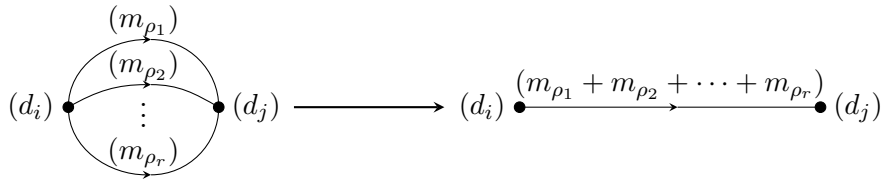
\begin{center}
    \begin{tikzpicture}[>=stealth]

      \filldraw (0,0) circle (2pt) (2,0) circle (2pt);

      \draw [->] [out=90, in=180] (0,0) to (1.1,0.9); \draw [out=0, in=90] (1.1,0.9) to (2,0);
      \draw [->] [out=30, in=180] (0,0) to (1.1,0.3); \draw [out=0, in=150] (1.1,0.3) to (2,0);
      \draw [->] [out=-90,in=180] (0,0) to (1.1,-0.9);\draw [out=0, in=-90] (1.1,-0.9) to (2,0);

      \begin{small}

	\draw (-0.5,0) node {$(d_i)$} (2.5,0) node{$(d_j)$};

	\draw (1,1.15) node {$(m_{\rho_1})$} (1,0.55) node {$(m_{\rho_2})$} (1,-0.6) node {$(m_{\rho_r})$};
	\draw (1,0) node {$\vdots$};

      \end{small}

      \draw [->, thick] (3,0) -- (5,0);

      \filldraw (6,0) circle (2pt) (10,0) circle (2pt);

      \draw[->] (6,0) -- (8.1,0); \draw (8.1,0) -- (10,0);

      \begin{small}

	\draw (5.5,0) node {$(d_i)$} (10.5,0) node{$(d_j)$};

	\draw (8,0.3) node{$(m_{\rho_1} + m_{\rho_2} + \dots + m_{\rho_r})$};

      \end{small}

    \end{tikzpicture}
  \end{center}
  The absolute valued quiver $\Gamma^C$ will be called the \emph{crushed (absolute valued) quiver}\index{Crushed!absolute valued quiver} of $\Gamma$.
\end{defn}

Note that $\Gamma^C$ does indeed satisfy the definition of an absolute valued quiver. Take any $\rho:i \rightarrow j \in \Gamma_1^C$. Then $d_i^C=d_i \; | \; m_\alpha$ for all $\alpha:i \rightarrow j$ in $\Gamma_1$. Thus, $d_i^C \; | \; \left( \sum_{(\alpha:i\rightarrow j) \; \in \; \Gamma_1} m_\alpha \right) = m_\rho^C$. The same is true for $d_j^C$. Therefore, $\Gamma^C$ is an absolute valued quiver.

The notion of crushing can be extended to relative valued quivers via the functor $\F$\index{$\F$}. Recall that if $\Gamma$ is an absolute valued quiver with valuation $(d_i,m_\rho)_{i \in \Gamma_0, \rho \in \Gamma_1}$, then $\F(\Gamma)$ is a relative valued quiver with valuation $(d_{ij}^\rho,d_{ji}^\rho)_{(\rho:i \rightarrow j) \in \F(\Gamma)_1}$ given by $d_{ij}^\rho = m_\rho / d_j$ and $d_{ji}^\rho = m_\rho / d_i$. So, the valuation of $\F(\Gamma^C)$ is given by
\[ (d^C)_{ij}^\rho = \left( \sum_{(\alpha:i \rightarrow j) \; \in \; \Gamma_1} m_\alpha \right) / d_j = \sum_{(\alpha:i \rightarrow j) \; \in \; \Gamma_1} d^\alpha_{ij} \] and likewise
\[ (d^C)_{ji}^\rho = \left( \sum_{(\alpha:i \rightarrow j) \; \in \; \Gamma_1} m_\alpha \right) / d_i = \sum_{(\alpha:i \rightarrow j) \; \in \; \Gamma_1} d^\alpha_{ji} \] for each $\rho:i \rightarrow j$ in $\F(\Gamma^C)_1$. We take this to be the definition of the crushed (relative valued) quiver of a relative valued quiver.

\begin{defn}[Crushed relative valued quiver] 
  Let $\Delta$ be a relative valued quiver\index{Quiver!valued!relative}. Define a new relative valued quiver, which we will denote $\Delta^C$, as follows:

  \begin{itemize}
    \item $\Delta_0^C = \Delta_0$,
    \item \# arrows $i \rightarrow j = 
            \begin{cases}
              1, & \text{if } \exists \; \rho:i \rightarrow j \in \Delta_1, \\
              0, & \text{otherwise},
            \end{cases}$
    \item $(d^C)^\rho_{ij} = \sum_{(\alpha:i \rightarrow j) \; \in \; \Delta_1} d^\alpha_{ij}$ and $(d^C)^\rho_{ji} = \sum_{(\alpha:i \rightarrow j) \; \in \; \Delta_1} d^\alpha_{ji}$ for all $\rho:i \rightarrow j$ in $\Delta_1^C$.
  \end{itemize}

  Again, the intuition is to ``crush'' all parallel arrows in $\Delta$ into a single arrow and sum the values.

  
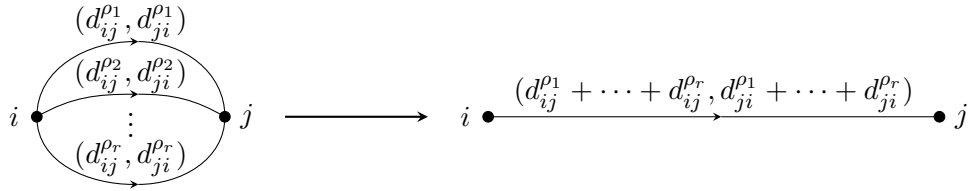
\begin{center}
    \begin{tikzpicture}[>=stealth]

      \filldraw (0,0) circle (2pt) (2.5,0) circle (2pt);

      \draw [->] [out=90, in=180] (0,0) to (1.35,1); \draw [out=0, in=90] (1.35,1) to (2.5,0);
      \draw [->] [out=30, in=180] (0,0) to (1.35,0.3); \draw [out=0, in=150] (1.35,0.3) to (2.5,0);
      \draw [->] [out=-90,in=180] (0,0) to (1.35,-0.9);\draw [out=0, in=-90] (1.35,-0.9) to (2.5,0);

      \begin{small};

	\draw (-0.3,0) node {$i$} (2.8,0) node{$j$};

	\draw (1.25,1.25) node {$(d^{\rho_1}_{ij},d^{\rho_1}_{ji})$} (1.25,0.55) node {$(d^{\rho_2}_{ij},d^{\rho_2}_{ji})$} (1.25,-0.55) node {$(d^{\rho_r}_{ij},d^{\rho_r}_{ji})$};
	\draw (1.25,0) node {$\vdots$};

      \end{small}

      \draw [->, thick] (3.3,0) -- (5.2,0);

      \filldraw (6,0) circle (2pt) (12,0) circle (2pt);

      \draw[->] (6,0) -- (9.1,0); \draw (9.1,0) -- (12,0);

      \begin{small}

	\draw (5.7,0) node {$i$} (12.3,0) node{$j$};

	\draw (9,0.3) node{$(d^{\rho_1}_{ij}+\dots+d^{\rho_r}_{ij},d^{\rho_1}_{ji}+\dots+d^{\rho_r}_{ji})$};

      \end{small}

    \end{tikzpicture}
  \end{center}
  The relative valued quiver $\Delta^C$ will be called the \emph{crushed (relative valued) quiver}\index{Crushed!relative valued quiver} of $\Delta$.
\end{defn}

\begin{defn}[Crushed species] 
  Let $\Q$ be a species\index{Species} with underlying relative valued quiver $\Delta$. Define a new species, which we will denote $\Q^C$, as follows:

  \begin{itemize}
    \item the underlying valued quiver of $\Q^C$ is $\Delta^C$,
    \item $K_i^C = K_i$ for all $i \in \Delta_0^C=\Delta_0$,
    \item $M_\rho^C = \bigoplus_{(\alpha:i\rightarrow j) \; \in \; \Delta_1} M_\alpha$ for all $\rho:i \rightarrow j$ in $\Delta_1^C$.
  \end{itemize}

  The intuition here is similar to that of the previous definitions; one ``crushes'' all bimodules along parallel arrows into a single bimodule by taking their direct sum. The species $\Q^C$ will be called the \emph{crushed species}\index{Crushed!species} of $\Q$.

\end{defn}

\begin{rmk}
  A \emph{crushed} $K$\emph{-species}\index{Crushed! $K$-species}\index{Species!$K$-species} is defined in exactly the same way, only using the crushed quiver of an absolute valued quiver instead of a relative valued quiver.
\end{rmk}

Note that in the above definition $\Q^C$ is indeed a species of $\Delta^C$. Clearly, $M_\rho^C$ is a $(K_j^C,K_i^C)$-bimodule for all $\rho:i \rightarrow j$ in $\Delta_1^C$. Moreover,
\begin{align*}
  \Hom_{K_j^C}(M_\rho^C,K_j^C) &= \Hom_{K_j}\left(\bigoplus_{(\alpha:i \rightarrow j) \in \Delta_1} M_\alpha, K_j\right) 
			       \cong \bigoplus_{(\alpha:i \rightarrow j) \in \Delta_1} \Hom_{K_j}(M_\alpha,K_j) \\
			       & \cong \bigoplus_{(\alpha:i \rightarrow j) \in \Delta_1} \Hom_{K_i}(M_\alpha,K_i) 
			       \cong \Hom_{K_i}\left(\bigoplus_{(\alpha:i \rightarrow j) \in \Delta_1} M_\alpha, K_i\right) \\
			       & = \Hom_{K_i}(M_\rho^C,K_i^C),
\end{align*}
where all isomorphisms are $(K_j^C,K_i^C)$-bimodule isomorphisms (in the second isomorphism we use the fact that $\Q$ is a species). Thus the duality condition for species holds. As for the dimension condition:
\begin{align*}
  \dim_{K_j^C}(M_\rho^C) &= \dim_{K_j}\left(\bigoplus_{(\alpha:i \rightarrow j) \in \Delta_1} M_\alpha \right) 
			 = \sum_{(\alpha:i \rightarrow j) \in \Delta_1} \dim_{K_j}(M_\alpha) 
			 = \sum_{(\alpha:i \rightarrow j) \in \Delta_1} d^\alpha_{ij} 
			 = (d^C)^\rho_{ij},
\end{align*}
where in the third equality we use the fact that $\Q$ is a species. Likewise, $\dim_{K_i^C}(M_\rho^C) = (d^C)^\rho_{ji}$. Thus, $\Q^C$ is a species of $\Delta^C$.

Note also that if $\Q$ is a $K$-species of an absolute valued quiver $\Gamma$, then $\Q^C$ is a $K$-species of $\Gamma^C$. Indeed, it is clear that $M_\rho^C$ is a $(K_j^C,K_i^C)$-bimodule on which $K$ acts centrally for all $\rho:i \rightarrow j$ in $\Gamma_1^C$ (since each summand satisfies this condition). Moreover, $\dim_K (K_i^C)=\dim_K (K_i)=d_i=d_i^C$ for all $i \in \Gamma_0^C$, thus the dimension criterion for the vertices is satisfied. Also, $\dim_K (M_\rho^C) = m_\rho^C$ by a computation similar to the above.

With the concept of crushed species/quivers, we obtain the following result, which gives a necessary and sufficient condition for two tensor rings/algebras\index{Tensor ring}\index{Tensor algebra} to be isomorphic.

\begin{theo}\label{theo:iso}
  Let $\Q$ and $\Q'$ be two species containing no oriented cycles. Then $T(\Q) \cong T(\Q')$\index{$T(\Q)$}\index{Tensor ring} as rings if and only if $\Q^C \cong \Q'^C$ as species. Moreover, if $\Q$ and $\Q'$ are $K$-species, then $T(\Q) \cong T(\Q)$\index{Tensor algebra} as $K$-algebras if and only if $\Q^C \cong \Q'^C$ as $K$-species.
\end{theo}

\begin{proof}
  Suppose $\Q$ and $\Q'$ are species with underlying relative valued quivers $\Delta$ and $\Delta'$, respectively. Throughout the proof we will use the familiar notation $D:= \Pi_{i \in \Delta_0} K_i$ and $M:=\bigoplus_{\rho \in \Delta_1} M_\rho$ (add primes for $\Q'$).

  The proof of the reverse implication is straightforward and so we leave the details to the reader.

  For the forward implication, assume $T(\Q) \cong T(\Q')$. Let $A=T(\Q)$ and $B=T(\Q')$ and let $\varphi:A \rightarrow B$ be a ring isomorphism. Then, there exists an induced ring isomorphism $\overline{\varphi}: A/ \rad A \rightarrow B/ \rad B$\index{Jacobson radical}\index{$\rad$}. By Lemma \ref{lem:rad}, $A / \rad A \cong D$ and $B / \rad B \cong D'$. Thus, we have an isomorphism $\widetilde{\varphi}_D : D \rightarrow D'$. It is not difficult to show that $\{1_{K_i}\}_{i \in \Delta_0}$ is the only complete set of primitive orthogonal idempotents in $D$ and that, likewise, $\{1_{K'_i}\}_{i \in \Delta'_0}$ is the only complete set of primitive orthogonal idempotents of $D'$. Since any isomorphism must bijectively map a complete set of primitive orthogonal idempotents to a complete set of primitive orthogonal idempotents, we may identify $\Delta_0$ and $\Delta'_0$ and assume, without loss of generality, that $\widetilde{\varphi}_D(1_{K_i}) = 1_{K'_i}$ for each $i \in \Delta_0=\Delta'_0$. Since $1_{K_i} \cdot D = K_i$ and $1_{K'_i} \cdot D' = K'_i$, we have that $\widetilde{\varphi}_D |_{K_i}$ is a ring isomorphism $K_i \rightarrow K'_i$ for each $i \in \Delta_0=\Delta'_0$.

  Now, $\varphi|_{\rad A}$ is a ring isomorphism from $\rad A$ to $\rad B$. Thus, as before, we have an induced ring isomorphism (and hence an abelian group isomorphism) $\overline{\varphi|_{\rad A}} : \rad A / (\rad A)^2 \rightarrow \rad B / (\rad B)^2$. Since $\rad A / (\rad A)^2 \cong M$ and $\rad B / (\rad B)^2 \cong M'$, we have an isomorphism $\widetilde{\varphi}_M : M \rightarrow M'$. For any $i,j \in \Delta_0$, $1_{K_j} \cdot M \cdot 1_{K_i} = \bigoplus_{(\rho:i \rightarrow j) \in \Delta_0} M_\rho=: {_jM_i}$ and $1_{K'_j} \cdot M' \cdot 1_{K'_i} = \bigoplus_{(\rho:i \rightarrow j) \in \Delta'_0} M'_\rho=: {_j M'_i}$. Therefore, $\widetilde{\varphi}_M |_{_j M_i}$ is an abelian group isomorphism ${_jM_i} \rightarrow {_jM'_i}$.

  Hence, $\{ \widetilde{\varphi}_D|_{K_i}, \widetilde{\varphi}_M|_{_jM_i} \}_{i \in \Delta_0^C=\Delta_0, (i \rightarrow j) \in \Delta_1^C}$ defines an isomorphism of species from $\Q^C$ to $\Q'^C$.

  In the case of $K$-species, one simply has to replace the terms ``ring'' with ``$K$-algebra'', ``ring morphism'' with ``$K$-algebra morphism'' and ``abelian group homomorphism'' with ``$K$-linear map'' and the proof is the same.
\end{proof}

If $\Q$ (and $\Q'$) contain oriented cycles, the arguments in the proof of Theorem \ref{theo:iso} fail since, in general, it is not true that $\rad T(\Q) = \bigoplus_{n=1}^\infty T^n(M)$. However, it seems likely that one could modify the proof to avoid using the radical. Hence, we offer the following conjecture.

\begin{conj}
  Theorem \ref{theo:iso} holds even if $\Q$ and $\Q'$ contain oriented cycles.
\end{conj}

\begin{rmk}\label{rmk:step1}
  Theorem \ref{theo:iso} serves as a first step in justifying Remark \ref{rmk:multi-arrows} (i.e.\ that we can always assume, without loss of generality, that we have no parallel arrows in our valued quivers) since a species with parallel arrows can always be crushed to one with only single arrows and its tensor algebra remains the same.
\end{rmk}

Theorem \ref{theo:iso} shows that there does not exist an equivalence on the level of valued quivers (relative or absolute) such that
\[ T(\Q) \cong T(\Q') \iff \Delta \text{ is equivalent to } \Delta' \]
since there are species (respectively $K$-species), with identical underlying valued quivers, that are not isomorphic as crushed species (respectively $K$-species) and hence have non-isomorphic tensor rings (respectively algebras) (see Example \ref{egg:kspecies}).

In the case of $K$-species, one may wonder what happens when we tensor $T(\Q)$ with the algebraic closure of $K$. Indeed, maybe we can find an equivalence on the level of absolute valued quivers such that
\[ \overline{K} \otimes_K T(\Q) \cong \overline{K} \otimes_K T(\Q') \iff \Gamma \text{ is equivalent to } \Gamma'. \]
The answer, unfortunately, is no. However, this idea does yield an interesting result. In \cite{hub}, Hubery showed that if $K$ is a finite field, then there is a field extension $F/K$ such that $F \otimes_K T(\Q)$ is isomorphic to the path algebra of a quiver. Our strategy of tensoring with the algebraic closure allows us to generalize this result for an arbitrary perfect field.

\begin{theo}\label{theo:tensor}
  Let $K$ be a perfect field\index{Perfect field} and $\Q$ be a $K$-species with underlying absolute valued quiver $\Gamma$ containing no oriented cycles such that $K_i$ is a field for each $i \in \Gamma_0$. Then $\overline{K} \otimes_K T(\Q)$\index{$T(\Q)$}\index{Tensor algebra} is isomorphic to the path algebra of a quiver\index{Path algebra}.
\end{theo}

\begin{proof}
  Let $A=\overline{K} \otimes_K T(\Q)$. Take any $i \in \Gamma_0$. It is a well-known fact that $\overline{K} \otimes_K K_i = \overline{K}^{d_i}$ (see for example \cite[Chapter V, Section 6, Proposition 2]{bourbaki} or the proof of \cite[Theorem 8.46]{jacobson}).

  Let $I = \overline{K} \otimes (\sum_{n=1}^\infty T^n(M))$. It is clear that $I$ is a two-sided ideal of $A$. Moreover, since $\Gamma$ has no cycles, $I$ is also nilpotent. Considering $A/I$, we see that
  \[ A/I \cong \overline{K} \otimes_K \left(\Pi_{i \in \Gamma_0} K_i \right) \cong \underbrace{\overline{K} \times \dots \times \overline{K}}_{\sum_{i \in \Gamma_0} d_i \text{ times}}. \]
  So, as in Lemma \ref{lem:rad}, $\rad A=I$\index{$\rad$}\index{Jacobson radical} and by \cite[Section I, Proposition 6.2]{ass}, $A$ is a basic\index{Basic algebra} finite-dimensional $\overline{K}$-algebra.

  We claim that $A$ is also hereditary\index{Hereditary algebra}. It is well-known that a ring is hereditary if and only if it is of global dimension at most 1. According to \cite[Theorem 16]{aus}, if $\Lambda_1$ and $\Lambda_2$ are $K$-algebras such that $\Lambda_1$ and $\Lambda_2$ are semiprimary (recall that a $K$-algebra $\Lambda$ is semiprimary if there is a two-sided nilpotent ideal $I$ such that $\Lambda/I$ is semisimple) and $(\Lambda_1/\rad \Lambda_1) \otimes_K (\Lambda_2 / \rad \Lambda_2)$ is semisimple, then $\gldim(\Lambda_1 \otimes_K \Lambda_2) = \gldim \Lambda_1 + \gldim \Lambda_2$. The $K$-algebras $\overline{K}$ and $T(\Q)$ satisfy these conditions. Indeed, $\overline{K}$ is simple and thus semiprimary. We know also that $T(\Q)/\rad T(\Q)$ is semisimple and $\rad T(\Q)$ is nilpotent since $\Gamma$ has no oriented cycles; thus $T(\Q)$ is semiprimary. Moreover,
  \begin{align*}
    (\overline{K} / \rad \overline{K}) \otimes_K (T(\Q)/ \rad T(\Q)) &\cong \overline{K} \otimes_K \left( \Pi_{i \in \Gamma_0} K_i \right) \\
    &\cong \underbrace{\overline{K} \times \cdots \times \overline{K}}_{\sum_{i \in \Gamma_0} d_i \text{ times}},
  \end{align*}
  which is semisimple.

  Therefore we have that:
  \[ \gldim A = \gldim(\overline{K} \otimes_K T(\Q)) = \gldim \overline{K} + \gldim T(\Q). \]
  However, $\gldim \overline{K} = 0$ (since all $\overline{K}$-modules are free) and $\gldim T(\Q) \le 1$ (since $T(\Q)$ is hereditary). Hence, $\gldim A \le 1$ and so $A$ is hereditary. By Theorem \ref{theo:pathalg}, $A$ is isomorphic to the path algebra of a quiver.
\end{proof}

\begin{rmk}
  Hubery goes further in \cite{hub}, constructing an automorphism $\sigma$ of the quiver $Q$ whose path algebra is isomorphic to $\overline{K}\otimes_K T(\Q)$ such that $(Q,\sigma)$ folds into $\Gamma$. It seems likely that this is possible here as well.
\end{rmk}

\begin{conj}
  Let $K$ be a perfect field, let $\Q$ be a $K$-species with underlying absolute valued quiver $\Gamma$ containing no oriented cycles such that $K_i$ is a field for each $i \in \Gamma_0$ and let $Q$ be a quiver such that $\overline{K} \otimes_K T(\Q) \cong \overline{K}Q$ (as in Theorem \ref{theo:tensor}). Then there exists an automorphism $\sigma$ of $Q$ such that $(Q,\sigma)$ folds into $\Gamma$.
\end{conj}

With Theorem \ref{theo:tensor}, we are able to use the methods of \cite[Chapter II, Section 3]{ass} to construct the quiver, $Q$, whose path algebra is isomorphic to $A=\overline{K} \otimes_K T(\Q)$. That is, the vertices of $Q$ are in one-to-one correspondence with $\{ \varepsilon_1, \dots, \varepsilon_n \}$, a complete set of primitive orthogonal idempotents\index{Idempotent} of $A$, and the number of arrows from the vertex corresponding to $\varepsilon_i$ to the vertex corresponding to $\varepsilon_j$ is given by $\dim_{\overline{K}} (\varepsilon_j \cdot (\rad A / \rad A^2) \cdot \varepsilon_i)$\index{Jacobson radical}\index{$\rad$}. We illustrate this in the next example.

\begin{egg}\label{egg:tensor}
  Let $\Gamma$ be the following absolute valued quiver.

  \begin{center}
    \begin{tikzpicture} [>=stealth]
      \draw (-1,0) node {$\Gamma$:};

      \filldraw (0,0) circle (2pt) (2,0) circle (2pt);

      \draw [->] (0,0) -- (1.1,0); \draw (1.1,0) -- (2,0);

      \begin{small}
	\draw (0,-0.3) node {(2)} (2,-0.3) node {(2)} (1,0.3) node {(4)};
      \end{small}
    \end{tikzpicture}
  \end{center}
  We can construct two $\mathbb{Q}$-species of $\Gamma$:

  \begin{center}
    \begin{tikzpicture} [>=stealth, baseline=0]
      \draw (-1,0) node {$\Q$:} (4,0) node {$\Q'$:};

      \filldraw (0,0) circle (2pt) (2,0) circle (2pt);

      \filldraw (5,0) circle (2pt) (7,0) circle (2pt);

      \draw [->] (0,0) -- (1.1,0); \draw (1.1,0) -- (2,0);

      \draw [->] (5,0) -- (6.1,0); \draw (6.1,0) -- (7,0);

      \begin{small}
	\draw (0,-0.3) node {$\mathbb{Q}(\sqrt{2})$} (1,0.3) node {$\mathbb{Q}(\sqrt{2})^2$} (2,-0.3) node {$\mathbb{Q}(\sqrt{2})$};

	\draw (5,-0.3) node {$\mathbb{Q}(\sqrt{2})$} (6,0.3) node {$\mathbb{Q}(\sqrt{2},\sqrt{3})$} (7,-0.3) node {$\mathbb{Q}(\sqrt{3})$};
      \end{small}
    \end{tikzpicture}.
  \end{center}
  Let $F=\overline{\mathbb{Q}}$, $A = F \otimes_\mathbb{Q} T(\Q)$ and $B = F \otimes_\mathbb{Q} T(\Q')$. We would like to find quivers $Q$ and $Q'$ with $A \cong FQ$ and $B \cong FQ'$.

  By direct computation, we see that $\varepsilon_1 = \frac{1}{2} ( (1\otimes1) + (\frac{1}{\sqrt{2}} \otimes \sqrt{2}) )$ and $\varepsilon_2 = \frac{1}{2} ( (1\otimes1) - (\frac{1}{\sqrt{2}} \otimes \sqrt{2}) )$ form a complete set of primitive orthogonal idempotents of $F \otimes_\mathbb{Q} \mathbb{Q}(\sqrt{2})$. Thus, $Q$ must have 4 vertices. To find the arrows, note that
  \begin{align*}
    \varepsilon_j \cdot (\rad A / \rad A^2) \cdot \varepsilon_i &= \varepsilon_j \cdot (F \otimes_\mathbb{Q} \mathbb{Q}(\sqrt{2})^2) \cdot \varepsilon_i \\
    &\cong \varepsilon_j \cdot (F \otimes_\mathbb{Q} \mathbb{Q}(\sqrt{2}))^2 \cdot \varepsilon_i \\
    &= F(\delta_{ij} \varepsilon_i, 0) \oplus F(0, \delta_{ij} \varepsilon_i),
  \end{align*}
  which has dimension $2$ if $i=j$ and $0$ otherwise. Hence, $A$ is isomorphic to the path algebra (over $F$) of

  \begin{center}
    \begin{tikzpicture} [>=stealth, baseline=-8]
      \draw (-1,-0.25) node {$Q:$};

      \filldraw (0,0) circle (2pt) (1,0) circle (2pt) (0,-0.5) circle (2pt) (1,-0.5) circle (2pt);

      \draw[->] (0,0.05) -- (0.6,0.05); \draw (0.6,0.05) -- (1,0.05);
      \draw[->] (0,-0.05) -- (0.6,-0.05); \draw (0.6,-0.05) -- (1,-0.05);

      \draw[->] (0,-0.45) -- (0.6,-0.45); \draw (0.6,-0.45) -- (1,-0.45);
      \draw[->] (0,-0.55) -- (0.6,-0.55); \draw (0.6,-0.55) -- (1,-0.55);
    \end{tikzpicture}$\;$.
  \end{center}

  For $Q'$, again $\{ \varepsilon_1, \varepsilon_2 \}$ is a complete set of primitive orthogonal idempotents of $F \otimes_\mathbb{Q} \mathbb{Q}(\sqrt{2})$. Likewise, $\zeta_1 = \frac{1}{2} ( (1\otimes1) + (\frac{1}{\sqrt{3}} \otimes \sqrt{3}) )$ and $\zeta_2 = \frac{1}{2} ( (1\otimes1) - (\frac{1}{\sqrt{3}} \otimes \sqrt{3}) )$ form a complete set of primitive orthogonal idempotents of $F \otimes_\mathbb{Q} \mathbb{Q}(\sqrt{3})$. So $Q'$ has 4 vertices. To find the arrows, note that
  \[ \zeta_j \cdot (\rad B / \rad B^2) \cdot \varepsilon_i = \zeta_j \cdot (F \otimes_\mathbb{Q} \mathbb{Q}(\sqrt{2},\sqrt{3})) \cdot \varepsilon_i \]
  and, since $\zeta_j \cdot \varepsilon_i \ne 0$ (as elements in $F \otimes_\mathbb{Q} \mathbb{Q}(\sqrt{2},\sqrt{3})$), this has dimension $1$ for all $i,j \in \{1,2\}$. Hence, $B$ is isomorphic to the path algebra (over $F$) of

  \begin{center}
    \begin{tikzpicture} [>=stealth, baseline=-8]
      \draw (-1,-0.25) node {$Q':$};

      \filldraw (0,0) circle (2pt) (1,0) circle (2pt) (0,-0.5) circle (2pt) (1,-0.5) circle (2pt);

      \draw[->] (0,0) -- (0.6,0); \draw (0.6,0) -- (1,0);
      \draw[->] (0,-0.5) -- (0.6,-0.5); \draw (0.6,-0.5) -- (1,-0.5);

      \draw[->] (0,0) -- (0.7,-0.35); \draw (0.7,-0.35) -- (1,-0.5);
      \draw[->] (0,-0.5) -- (0.7,-0.15); \draw (0.7,-0.15) -- (1,0);
    \end{tikzpicture}$\;$.
  \end{center}

\end{egg}

Notice that in Example \ref{egg:tensor}, $FQ$ and $FQ'$ are not isomorphic; this illustrates our earlier point; namely that there is no equivalence on the level of absolute valued quivers such that
\[ \overline{K} \otimes_K T(\Q) \cong \overline{K} \otimes_K T(\Q') \iff \Gamma \text{ is equivalent to } \Gamma'. \]
Therefore, it seems likely that Theorem \ref{theo:tensor} is the best that we can hope to achieve.

Note, however, that Theorem \ref{theo:tensor} fails if all the division rings in our $K$-modulation are not fields. Consider the following simple example.

\begin{egg}
  View $\mathbb{H}$, the quaternions, as an $\mathbb{R}$-species (that is, $\mathbb{H}$ is an $\mathbb{R}$-modulation of the absolute valued quiver with one vertex of value $4$ and no arrows). Consider the $\mathbb{C}$-algebra $\mathbb{C} \otimes_\mathbb{R} \mathbb{H}$. It is easy to see that $\mathbb{C} \otimes_\mathbb{R} \mathbb{H} \cong M_2(\mathbb{C})$, the algebra of $2$ by $2$ matrices with entries in $\mathbb{C}$. This algebra is not basic. Indeed, one can check (by direct computation) that
  \[ \left \{ \varepsilon_1 = \begin{pmatrix} 1 & 0 \\ 0 & 0 \end{pmatrix}, \varepsilon_2 = \begin{pmatrix} 0 & 0 \\ 0 & 1 \end{pmatrix} \right \} \]
  is a complete set of primitive orthogonal idempotents and that $M_2(\mathbb{C})\varepsilon_1 \cong M_2(\mathbb{C})\varepsilon_2 \cong \mathbb{C}^2$ as $M_2(\mathbb{C})$-modules. Thus, $\mathbb{C} \otimes_\mathbb{R} \mathbb{H}$ is not isomorphic to the path algebra of a quiver (since all path algebras are basic).
\end{egg}

While we cannot use Theorem \ref{theo:tensor} for arbitrary $K$-species, we do have the following.

\begin{cor}\label{cor:tensor}
  Let $K$ be a perfect field and $\Q$ be a $K$-species with underlying absolute valued quiver $\Gamma$ containing no oriented cycles. Then $\overline{K} \otimes_K T(\Q)$\index{$T(\Q)$}\index{Tensor algebra} is Morita equivalent\index{Morita equivalence} to the path algebra of a quiver\index{Path algebra}.
\end{cor}

\begin{proof}
  It suffices to show that $\overline{K} \otimes_K T(\Q)$ is hereditary (and then invoke Part 1 of Corollary \ref{cor:alg}). In the proof of Theorem \ref{theo:tensor}, all the arguments proving that $\overline{K} \otimes_K T(\Q)$ is hereditary go through as before, save for the proof that $(\overline{K}/ \rad \overline{K}) \otimes_K (T(\Q) / \rad T(\Q)) \cong \overline{K} \otimes_K (\Pi_{i \in \Gamma_0} K_i)$ is semisimple.

  To show this in the case that the $K_i$ are not necessarily all fields, pick some $i \in \Gamma_0$ and let $Z$ be the centre of $K_i$. Then
  \[ \overline{K} \otimes_K K_i \cong \overline{K} \otimes_K Z \otimes_Z K_i. \]
  The field $Z$ is a field extension of $K$ and so we may use the same arguments as in the proof of Theorem \ref{theo:tensor} to show $\overline{K} \otimes_K Z \cong \overline{K} \times \dots \times \overline{K}$. So
  \begin{align*}
    \overline{K} \otimes_K K_i &\cong (\overline{K} \times \dots \times \overline{K}) \otimes_Z K_i \\
		    &\cong (\overline{K} \otimes_Z K_i) \times \dots \times (\overline{K} \otimes_Z K_i).
  \end{align*}
  One can show that $\overline{K} \otimes_Z K_i \cong M_n(\overline{K})$ for some $n$, which is a simple ring. Thus, $\overline{K} \otimes_K K_i$ is semisimple, meaning that $\overline{K} \otimes_K (\Pi_{i \in \Gamma_0} K_i)$ is semisimple, completing the proof.
\end{proof}

\section{Representations}
In this section, we begin by defining representations of quivers and species. We will then see (Proposition \ref{prop:equ}) that representations of species (resp.\ quivers) are equivalent to modules of the corresponding tensor ring (resp.\ path algebra). This fact together with Section 3 (specifically Theorems \ref{theo:pathalg} and \ref{theo:tensalg}, and Corollary \ref{cor:alg}) shows why representations of quivers/species are worth studying; they allow us to understand the representations of any finite-dimensional algebra over a perfect field. We then discuss the root system associated to a valued quiver, which encodes a surprisingly large amount of information about the representation theory of species (see Theorems \ref{theo:main} and \ref{theo:kac}, and Proposition \ref{prop:dx}). From Section 1, we know that every valued quiver can be obtained by folding a quiver with automorphism. Thus, we end the section with a discussion on how much of the data of the representation theory of a species is contained in a corresponding quiver with automorphism.

Throughout this section, we make the assumption (unless otherwise specified) that all quivers/species are connected and contain no oriented cycles. Also, whenever there is no need to distinguish between relative or absolute valued quivers, we will simply use the term ``valued quiver''\index{Quiver!valued} and denote it by $\Omega$. We let $\{e_i\}_{i \in \Omega_0}$ be the standard basis of $\mathbb{Z}^{\Omega_0}$ for a valued quiver $\Omega$.

\begin{defn}[Representation of a quiver] 
  A \emph{representation}\index{Representation!of a quiver} $V=(V_i,f_\rho)_{i\in Q_0, \rho \in Q_1}$ of a quiver $Q$ over the field $K$ consists of a $K$-vector space $V_i$ for each $i \in Q_0$ and a $K$-linear map
  \[ f_\rho:V_{t(\rho)} \rightarrow V_{h(\rho)}, \]
  for each $\rho \in Q_1$. If each $V_i$ is finite-dimensional, we call $\dimbar V = (\dim_K V_i)_{i \in Q_0} \in \mathbb{N}^{Q_0}$\index{$\dimbar$} the \emph{graded dimension}\index{Graded dimension} of $V$.

  A \emph{morphism} of $Q$ representations 
  \[ \varphi: V=(V_i,f_\rho)_{i\in Q_0, \rho \in Q_1} \rightarrow W=(W_i,g_\rho)_{i\in Q_0, \rho \in Q_1} \]
  consists of a $K$-linear map $\varphi_i : V_i \rightarrow W_i$ for each $i \in Q_0$ such that $\varphi_{h(\rho)} \circ f_\rho = g_\rho \circ \varphi_{t(\rho)}$ for all $\rho \in Q_1$. That is, the following diagram commutes for all $\rho \in Q_1$.

  \begin{center}
    \begin{tikzpicture}
      \draw (0,0) node {$V_{t(\rho)}$};
      \draw[->] (0.4,0) -- (1.4,0); \draw (0.9,0.3) node {$f_\rho$};
      \draw (1.8,0) node {$V_{h(\rho)}$};
      \draw[->] (0,-0.3) -- (0,-1.3); \draw (-0.4,-0.8) node {$\varphi_{t(\rho)}$};
      \draw (0,-1.6) node {$W_{t(\rho)}$};
      \draw[->] (0.4,-1.6) -- (1.4,-1.6); \draw (0.9,-1.9) node {$g_\rho$};
      \draw (1.9,-1.6) node {$W_{h(\rho)}$};
      \draw[->] (1.9,-0.3) -- (1.9,-1.3); \draw (2.4,-0.8) node {$\varphi_{h(\rho)}$};
    \end{tikzpicture}
  \end{center}
  We let $\R_K(Q)$\index{$\R_K(Q)$} denote the category of finite-dimensional representations of $Q$ over $K$.
\end{defn}

\begin{defn}[Representation of a species] 
  A \emph{representation}\index{Representation!of a species} $V=(V_i,f_\rho)_{i\in \Delta_0, \rho \in \Delta_1}$ of a species (or $K$-species) $\Q$ consists of a $K_i$-vector space $V_i$ for each $i \in \Delta_0$ and a $K_{h(\rho)}$-linear map 
  \[ f_\rho:M_\rho \otimes_{K_{t(\rho)}} V_{t(\rho)} \rightarrow V_{h(\rho)}, \]
  for each $\rho \in \Delta_1$. If all $V_i$ are finite-dimensional (over their respective rings), we call $\dimbar V = (\dim_{K_i} V_i)_{i \in \Delta_0} \in \mathbb{N}^{\Delta_0}$\index{$\dimbar$} the \emph{graded dimension}\index{Graded dimension} of $V$.

  A \emph{morphism} of $\Q$ representations 
  \[ \varphi : V=(V_i,f_\rho)_{i\in \Delta_0, \rho \in \Delta_1} \rightarrow W=(W_i,g_\rho)_{i\in \Delta_0, \rho \in \Delta_1} \]
  consists of a $K_i$-linear map $\varphi_i: V_i \rightarrow W_i$ for each $i \in \Delta_0$ such that $\varphi_{h(\rho)} \circ f_\rho = g_\rho \circ (\id_{M_\rho} \otimes \varphi_{t(\rho)})$ for all $\rho \in \Delta_1$. That is, the following diagram commutes for all $\rho \in \Delta_1$.

  \begin{center}
    \begin{tikzpicture}
      \draw (0,0) node {$M_\rho \otimes_{K_{t(\rho)}} V_{t(\rho)}$};
      \draw[->] (1.3,0) -- (2.3,0); \draw (1.8,0.3) node {$f_\rho$};
      \draw (2.8,0) node {$V_{h(\rho)}$};
      \draw[->] (0,-0.3) -- (0,-1.3); \draw (-1.1,-0.8) node {$\id_{M_\rho} \otimes \varphi_{t(\rho)}$};
      \draw (0,-1.6) node {$M_\rho \otimes_{K_{t(\rho)}} W_{t(\rho)}$};
      \draw[->] (1.4,-1.6) -- (2.2,-1.6); \draw (1.8,-1.9) node {$g_\rho$};
      \draw (2.8,-1.6) node {$W_{h(\rho)}$};
      \draw[->] (2.8,-0.3) -- (2.8,-1.3); \draw (3.4,-0.8) node {$\varphi_{h(\rho)}$};
    \end{tikzpicture}
  \end{center}
  We let $\R(\Q)$\index{$\R(\Q)$} denote the category of finite-dimensional representations of $\Q$. If $\Q$ is a $K$-species, we use the notation $\R_K(\Q)$\index{$\R_K(\Q)$}.
\end{defn}

Note that if $\Q$ is a $K$-species\index{Species!$K$-species} of a trivially valued absolute valued quiver $\Gamma$, then, as before, all $K_i \cong K$ (as $K$-algebras) and all $M_\rho \cong K$ (as bimodules). Thus, a representation of $\Q$ is a representation of the underlying (non-valued) quiver of $\Gamma$. Therefore, by viewing quivers as trivially valued absolute valued quivers, representations of species are a generalization of representations of quivers.

It is well-known that, for a quiver $Q$, the category $\R_K(Q)$\index{$\R_K(Q)$} is equivalent to $KQ$-mod\index{$KQ$}\index{Path algebra}, the category of finitely-generated (left) $KQ$-modules. This fact generalizes nicely for species.

\begin{prop}\label{prop:equ}
  Let $\Q$ be a species (possibly with oriented cycles). Then $\R(\Q)$\index{$\R(\Q)$} is equivalent to $T(\Q)$-$\mathrm{mod}$\index{$T(\Q)$}\index{Tensor ring}.
\end{prop}

\begin{proof}
  See \cite[Proposition 10.1]{dr2}. While the proof there is given only for $K$-species, the same arguments hold for species in general.
\end{proof}

\begin{rmk}
  Proposition \ref{prop:equ}, together with Theorem \ref{theo:iso}, justifies Remark \ref{rmk:multi-arrows} (i.e.\ that we can always assume, without loss of generality, that our valued quivers contain no parallel arrows) since a species with parallel arrows can always be crushed\index{Crushed!species} to one with only single arrows and its tensor algebra remains the same. Since $T(\Q)$-mod is equivalent to $\R(\Q)$-mod, the representation theory of any species is equivalent to the representation of a species with only single arrows (its crushed species). While allowing parallel arrows in our definition of species is not necessary, there are situations where it may be advantageous as the next example demonstrates.
\end{rmk}

\begin{egg}\label{egg:just}
  Let $\Delta$ be the following valued quiver.

  
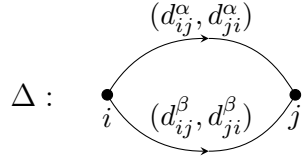
\begin{center}
    \begin{tikzpicture}[>=stealth]

      \draw (-1,0) node {$\Delta:$};

      \filldraw (0,0) circle (2pt) (2.5,0) circle (2pt);

      \draw [->] [out=60, in=180] (0,0) to (1.35,0.75); \draw [out=0, in=120] (1.35,0.75) to (2.5,0);
      \draw [->] [out=-60,in=180] (0,0) to (1.35,-0.75);\draw [out=0, in=-120] (1.35,-0.75) to (2.5,0);

      \begin{small}

	\draw (0,-0.3) node {$i$} (2.5,-0.3) node{$j$};

	\draw (1.25,1) node {$(d^{\alpha}_{ij},d^{\alpha}_{ji})$} (1.25,-0.3) node {$(d^{\beta}_{ij},d^{\beta}_{ji})$};

      \end{small}
    \end{tikzpicture}
  \end{center}
  Then

  
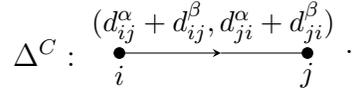
\begin{center}
    \begin{tikzpicture}[>=stealth, baseline=0]

      \draw (-1,0) node {$\Delta^C:$};

      \filldraw (0,0) circle (2pt) (2.5,0) circle (2pt);

      \draw [->] (0,0) to (1.35,0); \draw (1.35,0) to (2.5,0);

      \begin{small}

	\draw (0,-0.3) node {$i$} (2.5,-0.3) node{$j$};

	\draw (1.25,0.4) node {$(d^{\alpha}_{ij}+d^{\beta}_{ij},d^{\alpha}_{ji}+d^{\beta}_{ji})$};

      \end{small}
    \end{tikzpicture}.
  \end{center}
  Any modulation of $\Delta$,

  
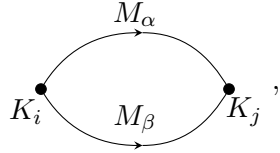
\begin{center}
    \begin{tikzpicture}[>=stealth, baseline=0]

      \filldraw (0,0) circle (2pt) (2.5,0) circle (2pt);

      \draw [->] [out=60, in=180] (0,0) to (1.35,0.75); \draw [out=0, in=120] (1.35,0.75) to (2.5,0);
      \draw [->] [out=-60,in=180] (0,0) to (1.35,-0.75);\draw [out=0, in=-120] (1.35,-0.75) to (2.5,0);

      \begin{small}

	\draw (-0.2,-0.3) node {$K_i$} (2.7,-0.3) node{$K_j$};

	\draw (1.25,1) node {$M_\alpha$} (1.25,-0.4) node {$M_\beta$};

      \end{small}
    \end{tikzpicture},
  \end{center}
  yields a modulation of $\Delta^C$,

  
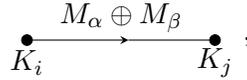
\begin{center}
    \begin{tikzpicture}[>=stealth, baseline=0]

      \filldraw (0,0) circle (2pt) (2.5,0) circle (2pt);

      \draw [->] (0,0) to (1.35,0); \draw (1.35,0) to (2.5,0);

      \begin{small}

	\draw (0,-0.3) node {$K_i$} (2.5,-0.3) node{$K_j$};

	\draw (1.25,0.3) node {$M_\alpha \oplus M_\beta$};

      \end{small}
    \end{tikzpicture},
  \end{center}
  and the representation theory of both these species is identical. However, the converse is not true. That is, not every modulation of $\Delta^C$ yields a modulation of $\Delta$. For example, one can choose a modulation

  
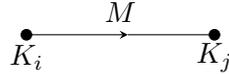
\begin{center}
    \begin{tikzpicture}[>=stealth]

      \filldraw (0,0) circle (2pt) (2.5,0) circle (2pt);

      \draw [->] (0,0) to (1.35,0); \draw (1.35,0) to (2.5,0);

      \begin{small}

	\draw (0,-0.3) node {$K_i$} (2.5,-0.3) node{$K_j$};

	\draw (1.25,0.3) node {$M$};

      \end{small}
    \end{tikzpicture}
  \end{center}
  such that $M$ is indecomposable, and thus cannot be written as $M=M_1 \oplus M_2$ (with $M_1, M_2 \ne 0$) to yield a modulation of $\Delta$. Thus, we can think of modulations of $\Delta$ as being ``special'' modulations of $\Delta^C$ where the bimodule attached to its arrow can be written (nontrivially) as the direct sum of two bimodules.
\end{egg}

Example \ref{egg:just} illustrates why one may wish to allow parallel arrows in the definition of species; they may be used as a way of ensuring that the bimodules in our modulation decompose into a direct sum of proper sub-bimodules.

\begin{defn}[Indecomposable representation] 
  Let $V=(V_i,f_\rho)_{i\in \Delta_0, \rho \in \Delta_1}$ and \linebreak $W=(W_i,g_\rho)_{i\in \Delta_0, \rho \in \Delta_1}$ be representations of a species (or a quiver). The \emph{direct sum}\index{Representation!direct sum} of $V$ and $W$ is
  \[ V \oplus W = (V_i \oplus W_i, f_\rho \oplus g_\rho)_{i \in \Delta_0, \rho \in \Delta_1}. \]
  A representation $U$ is said to be \emph{indecomposable}\index{Representation!indecomposable} if $U=V \oplus W$ implies $V=U$ or $W=U$.
\end{defn}

Because we restrict ourselves to finite-dimensional representations, the Krull-Schmidt theorem holds. That is, every representation can be written uniquely as a direct sum of indecomposable representations (up to isomorphism and permutation of the components). Thus, the study of all representations of a species (or quiver) reduces to the study of its indecomposable representations.

We say that a species/quiver is of \emph{finite} representation type\index{Representation type!finite} if it has only finitely many non-isomorphic indecomposable representations. It is of \emph{tame}\index{Representation type!tame}\index{Representation type!affine} (or \emph{affine}) representation type if it has infinitely many non-isomorphic indecomposable representations, but they can be divided into finitely many one parameter families. Otherwise, it is of \emph{wild}\index{Representation type!wild} representation type.

Thus the natural question to ask is: can we classify all species/quivers of finite type, tame type and wild type? The answer, as it turns out, is yes. However, we first need a few additional concepts.

\begin{defn}[Euler, symmetric Euler and Tits forms]\label{defn:euler} 
  The \emph{Euler form}\index{Euler form} of an absolute valued quiver\index{Quiver!valued!absolute} $\Gamma$ with valuation $(d_i,m_\rho)_{i \in \Gamma_0, \rho \in \Gamma_1}$ is the bilinear form $\langle-,-\rangle:\mathbb{Z}^{\Gamma_0} \times \mathbb{Z}^{\Gamma_0} \rightarrow \mathbb{Z}$ given by:
  \[ \langle x,y \rangle=\sum_{i\in \Gamma_0} d_i x_i y_i -\sum_{\rho \in \Gamma_1} m_\rho x_{t(\rho)} y_{h(\rho)}. \]
  The \emph{symmetric Euler form}\index{Euler form!symmetric} $(-,-):\mathbb{Z}^{\Gamma_0} \times \mathbb{Z}^{\Gamma_0} \rightarrow \mathbb{Z}$ is given by:
  \[ (x,y)=\langle x,y \rangle +\langle y,x \rangle. \]
  The \emph{Tits form}\index{Tits form} $q:\mathbb{Z}^{\Gamma_0} \rightarrow \mathbb{Z}$ is given by:
  \[ q(x)=\langle x,x \rangle. \]
\end{defn}

\begin{rmk}
  If we take $\Gamma$ to be trivially valued (i.e.\ all $d_i=m_\rho=1$), we recover the usual definitions of these forms for quivers (see, for example \cite[Definitions 3.6.7, 3.6.8, 3.6.9]{cnp}).
\end{rmk}

\begin{rmk}
  Notice that the symmetric Euler form and the Tits form do not depend on the orientation of our quiver.
\end{rmk}

\begin{rmk}
  Given a relative valued quiver $\Delta$\index{Quiver!valued!relative}, we have seen in Lemma \ref{lem:functorF} that we can choose an absolute valued quiver $\Gamma$ such that $\F(\Gamma)=\Delta$\index{$\F$} (this is equivalent to making a choice of positive integers $f_i$ in Definition \ref{defn:rvq}) with $d^\rho_{ij}=m_\rho / d_j$ and $d^\rho_{ji}=m_\rho / d_i$ for all $\rho:i \rightarrow j$ in $\Delta_1$. It is easy to see that (as long as the quiver is connected) for any other absolute valued quiver\index{Quiver!valued!absolute} $\Gamma'$ with $\F(\Gamma')=\Delta$, there is a $\lambda \in \mathbb{Q}^+$ such that $d'_i = \lambda d_i$ for all $i \in \Delta_0$. Thus, we define the Euler, symmetric Euler and Tits forms on $\Delta$ to be the corresponding forms on $\Gamma$, which are well-defined up to positive rational multiple.
\end{rmk}

\begin{defn}[Generalized Cartan matrix]\label{defn:cartan}
  Let $\mathcal{I}$ be an indexing set. A \emph{generalized Cartan matrix}\index{Generalized Cartan matrix} $C=(c_{ij})$, $i,j \in \mathcal{I}$, is an integer matrix satisfying:

    \begin{itemize}
      \item $c_{ii} = 2$, for all $i \in \mathcal{I}$;

      \item $c_{ij} \le 0$, for all $i \ne j \in \mathcal{I}$;

      \item $c_{ij}=0 \iff c_{ji}=0$, for all $i,j \in \mathcal{I}$.
    \end{itemize}

  A generalized Cartan matrix $C$ is \emph{symmetrizable} if there exists a diagonal matrix $D$ (called the \emph{symmetrizer}) such that $DC$ is symmetric.
\end{defn}

Note that, for any valued quiver\index{Quiver!valued} $\Omega$, $c_{ij}=2\dfrac{(e_i,e_j)}{(e_i,e_i)}$ defines a generalized Cartan matrix, since $(e_i,e_i)=2d_i$ and $(e_i,e_j)=-\sum_\rho m_\rho$ for $i\ne j$, where the sum is taken over all arrows between $i$ and $j$ (regardless of orientation). So,
\[ c_{ij}=\begin{cases} 2, & \text{if }i=j, \\ -\sum_\rho m_\rho / d_i = -\sum_\rho d_{ij}^\rho, & \text{if }i \ne j. \end{cases}\]
From this we see that two valued quivers $\Omega$ and $\Omega'$ have the same generalized Cartan matrix (up to ordering of the rows and columns) if and only if $\Omega^C \cong \Omega'^C$ as relative valued quivers (by this we mean that if $\Omega$ and $\Omega'$ are relative valued quivers, then $\Omega^C \cong \Omega'^C$ and if they are absolute valued quivers, then $\F(\Omega)^C \cong \F(\Omega')^C$\index{$\F$}). If all $d_i$ are equal (or alternatively, $\sum_\rho d_{ij}^\rho = \sum_\rho d_{ji}^\rho$ for all adjacent $i$ and $j$), then the matrix is symmetric, otherwise it is symmetrizable with symmetrizer $D=\diag(d_i)_{i \in \Omega_0}$. Moreover, every symmetrizable Cartan matrix can be obtained in this way. This is one of the motivations for working with species. When working with species we can obtain non-symmetric Cartan matrices, but when restricted to quivers, only symmetric Cartan matrices arise. For every generalized Cartan matrix, we have its associated Kac-Moody Lie algebra.

\begin{defn}[Kac-Moody Lie algebra]\label{defn:kac-moody}
  Let $C=(c_{ij})$ be an $n \times n$ generalized Cartan matrix. Then the \emph{Kac-Moody Lie algebra}\index{Kac-Moody Lie algebra} of $C$ is the complex Lie algebra generated by $e_i,f_i,h_i$ for $1 \le i \le n$, subject to the following relations.
  \begin{itemize}
    \item $[h_i,h_j]=0$ for all $i,j$,
    \item $[h_i,e_j]=c_{ij}e_j$ and $[h_i,f_j]=-c_{ij}f_j$ for all $i,j$,
    \item $[e_i,f_i]=h_i$ for each $i$ and $[e_i,f_j]=0$ for all $i \ne j$,
    \item $(\ad e_i)^{1-c_{ij}}(e_j)=0$ and $(\ad f_i)^{1-c_{ij}}(f_j)$ for all $i \ne j$.
  \end{itemize}
\end{defn}

Therefore, to every valued quiver, we can associate a generalized Cartan matrix and its corresponding Kac-Moody Lie algebra. It is only fitting then, that we discuss root systems.

\begin{defn}[Root system of a valued quiver]

  Let $\Omega$ be a valued quiver\index{Quiver!valued}.

  \begin{itemize}

    \item For each $i \in \Omega_0$, define the \emph{simple reflection}\index{Simple reflection} through $i$ to be the linear transformation $r_i:\mathbb{Z}^{\Omega_0} \rightarrow \mathbb{Z}^{\Omega_0}$ given by:
    \[ r_i(x)=x-2\dfrac{(x,e_i)}{(e_i,e_i)}e_i. \]
                                                      
    \item The \emph{Weyl group}\index{Weyl group}, which we denote by $\mathcal{W}$, is the subgroup of $\Aut(\mathbb{Z}^{\Omega_0})$ generated by the simple reflections $r_i$, $i\in 	\Omega_0$. 

    \item An element $x \in \mathbb{Z}^{\Omega_0}$ is called a \emph{real root}\index{Root!real} if $\exists \; w \in \mathcal{W}$ such that $x=w(e_i)$ for some $i \in \Omega_0$.

    \item The \emph{support}\index{Support} of an element $x \in \mathbb{Z}^{\Omega_0}$ is defined as $\supp(x)=\{i \in \Omega_0 \; | \; x_i \ne 0\}$ and we say $\supp(x)$ is connected if the full subquiver of $\Omega$ with vertex set $\supp(x)$ is connected. Then the \emph{fundamental set}\index{Fundamental set} is defined as $\mathcal{F}=\{0 \ne x \in \mathbb{N}^{\Omega_0} \; | \; (x,e_i) \le 0 \text{ for all } i \in \Omega_0 \text{ and } \supp(x) \text{ is connected}\}$.

    \item An element $x \in \mathbb{Z}^{\Omega_0}$ is called an \emph{imaginary root}\index{Root!imaginary} if $x \in \bigcup_{w \in \mathcal{W}} w(\mathcal{F}) \cup w(-\mathcal{F})$.

    \item The \emph{root system}\index{Root!system} of $\Omega$, denoted $\Phi(\Omega)$\index{$\Phi(\Omega)$} is the set of all real and imaginary roots.

    \item We call a root $x$ \emph{positive}\index{Root!positive} (resp. \emph{negative}\index{Root!negative}) if $x_i \ge 0$ (resp. $x_i \le 0$) $\forall \; i \in \Omega_0$. We write $\Phi^+(\Omega)$ for the set of positive roots and $\Phi^-(\Omega)$ for the set of negative roots.

  \end{itemize}
\end{defn}

\begin{defn}[Stable element]
  An element $x \in \mathbb{Z}^{\Omega_0}$ is called \emph{stable}\index{Stable element} if $w(x)=x$ for all $w \in \mathcal{W}$.
\end{defn}

\begin{rmk}
  It is worth noting that, while a stable element need not be an imaginary root (see \cite[Example 6.15]{me}), it is always the sum of imaginary roots (see \cite[Lemma 6.16]{me}).
\end{rmk}

\begin{defn}[Discrete and continuous dimension types]
  An indecomposable representation\index{Representation!indecomposable} $V$ of a species (or quiver) is of \emph{discrete dimension type}\index{Dimension type!discrete} if it is the unique indecomposable representation (up to isomorphism) with graded dimension $\dimbar V$. Otherwise, it is of \emph{continuous dimension type}\index{Dimension type!continuous}.
\end{defn}

With all these concepts in mind, we can neatly classify all species of finite and tame representation type. Note that in the case of quivers, this was originally done by Gabriel (see \cite{g}). It was later generalized to species by Dlab and Ringel.

\begin{theo}\textbf{\emph{\cite[Main Theorem]{dr}}}\label{theo:main}
  Let $\mathcal{Q}$ be species\index{Species} of a connected relative valued quiver\index{Quiver!valued!relative} $\Delta$. Then:

  \begin{enumerate}
    \item$\mathcal{Q}$ is of finite representation type\index{Representation type!finite} if and only if the underlying undirected valued graph of $\Delta$ is a Dynkin diagram\index{Dynkin diagram} of finite type (see \cite{dr} for a list of the Dynkin diagrams). Moreover, $\dimbar:\R(\Q) \rightarrow \mathbb{Z}^{\Delta_0}$\index{Graded dimension}\index{$\dimbar$} induces a bijection between the isomorphism classes of the indecomposable representations\index{Representation!indecomposable} of $\Q$ and the positive real roots\index{Root!positive}\index{Root!real} of its root system\index{Root!system}.

    \item If the underlying undirected valued graph of $\Delta$ is an extended Dynkin diagram\index{Dynkin diagram!extended} (see \cite{dr} for a list of the extended Dynkin diagrams), then $\dimbar: \R(\Q) \rightarrow \mathbb{Z}^{\Delta_0}$\index{Graded dimension}\index{$\dimbar$} induces a bijection between the isomorphism classes of the indecomposable representations\index{Representation!indecomposable} of $\Q$ of discrete dimension type\index{Dimension type!discrete} and the positive real roots\index{Root!positive}\index{Root!real} of its root system\index{Root!system}. Moreover, there exists a unique stable element\index{Stable element} (up to rational multiple) $n \in \Phi(\Q)$ and the indecomposable representations of continuous dimension type\index{Dimension type!continuous} are those whose graded dimension is a positive multiple of $n$. If $\Q$ is a $K$-species\index{Species!$K$-species}, then $\Q$ is of tame representation type\index{Representation type!tame} if and only if the underlying undirected valued graph of $\Delta$ is an extended Dynkin diagram.
  \end{enumerate}
\end{theo}

\begin{rmk}
  See \cite[p.\ 57]{dr} (and \cite{r2}) for a proof that a $K$-species $\Q$ is tame if and only if $\Delta$ is an extended diagram.
\end{rmk}

\begin{rmk}
  In the case that the underlying undirected valued graph of $\Delta$ is an extended Dynkin diagram, the indecomposable representations of continuous dimension type of $\Q$ can be derived from the indecomposable representations of continuous dimension type of a suitable species with underlying undirected valued graph $\mathrm{\widetilde{A}}_{11}$ or $\mathrm{\widetilde{A}}_{12}$ (see \cite[Theorem 5.1]{dr}).
\end{rmk}

Theorem \ref{theo:main} shows a remarkable connection between the representation theory of species and the theory of root systems of Lie algebras. In the case of quivers, Kac was able to show that this connection is stronger still.

\begin{theo}\textbf{\emph{\cite[Theorems 2 and 3]{kac2}}}\label{theo:kac}
  Let $Q$ be a quiver\index{Quiver} with no loops (though possibly with oriented cycles) and $K$ an algebraically closed field. Then there is an indecomposable representation\index{Representation!indecomposable} of $Q$ of graded dimension\index{Graded dimension} $\alpha$ if and only if $\alpha \in \Phi^+(Q)$. Moreover, if $\alpha$ is a real positive root\index{Root!real}\index{Root!positive}, then there is a unique indecomposable representation of $Q$ (up to isomorphism) of graded dimension $\alpha$. If $\alpha$ is an imaginary positive root\index{Root!imaginary}\index{Root!positive}, then there are infinitely many non-isomorphic indecomposable representations of $Q$ of graded dimension $\alpha$.
\end{theo}

It is not known whether Kac's theorem generalizes fully for species\index{Species}, however, it does for certain classes of species. Indeed, in the case of a species of finite or tame representation type, one can apply Theorem \ref{theo:main}. In the case of $K$-species when $K$ is a finite field, we have the following result by Deng and Xiao.

\begin{prop}\textbf{\emph{\cite[Proposition 3.3]{dx1}}}\label{prop:dx}
  Let $\Q$ be a $K$-species\index{Species!$K$-species} ($K$ a finite field) containing no oriented cycles. Then there exists an indecomposable representation\index{Representation!indecomposable} of $\Q$ of graded dimension $\alpha$ if and only if $\alpha \in \Phi^+(\Q)$. Moreover, if $\alpha$ is a real positive root\index{Root!positive}\index{Root!real}, then there is a unique indecomposable representation of $\Q$ (up to isomorphism) of graded dimension $\alpha$.
\end{prop}

Based on these results, we see that much of the information about the representation theory of a species is encoded in its underlying valued quiver/graph. Recall from Section 1 that any valued quiver can be obtained by folding\index{Folding} a quiver with automorphism. So, one may ask: how much information is encoded in this quiver with automorphism?

We continue our assumption that $Q$ contains no oriented cycles; however for what follows this is more restrictive than we need. It would be enough to assume that $Q$ contains no loops and that no arrow connects two vertices in the same $\sigma$-orbit (see \cite[Lemma 6.24]{me} for a proof that this is indeed a weaker condition).

Suppose $(Q,\sigma)$ is a quiver with automorphism and let $V=(V_i,f_\rho)_{i \in Q_0, \rho \in Q_1}$ be a representation of $Q$. Define a new representation $V^\sigma = (V^\sigma_i,f^\sigma_i)_{i \in Q_0, \rho \in Q_1}$ by $V^\sigma_i = V_{\sigma^{-1}(i)}$ and $f^\sigma_\rho = f_{\sigma^{-1}(\rho)}$.

\begin{defn}[Isomorphically invariant representation]
  Let $(Q,\sigma)$ be a quiver\index{Quiver} with automorphism. A representation $V=(V_i,f_\rho)_{i \in Q_0, \rho \in Q_1}$ is called \emph{isomorphically invariant}\index{Representation!(isomorphically) invariant} (or simply \emph{invariant}) if $V^\sigma \cong V$ as representations of $Q$.

  We say an invariant representation $V$ is \emph{invariant-indecomposable}\index{Representation!invariant-indecomposable} if $V = W_1 \oplus W_2$ such that $W_1$ and $W_2$ are invariant representations implies $W_1=V$ or $W_2=V$.
\end{defn}

It is not hard to see that the invariant-indecomposable representations are precisely those of the form
\[ V = W \oplus W^\sigma \oplus \dots \oplus W^{\sigma^{r-1}} \]
where $W$ is an indecomposable representation and $r$ is the least positive integer such that $W^{\sigma^r} \cong W$.

Let $(\mathbb{Z}^{Q_0})^\sigma = \{ \alpha \in \mathbb{Z}^{Q_0} \; | \; \alpha_i = \alpha_j \text{ for all } i \text{ and } j \text{ in the same orbit}\}$. Suppose $(Q,\sigma)$ folds into $\Omega$ and write $\overline{i} \in \Omega_0$ for the orbit of $i \in Q_0$. We then have a well-defined function\index{$\f$}
\[ \f : (\mathbb{Z}^{Q_0})^\sigma \rightarrow \mathbb{Z}^{\Omega_0} \]
defined by $\f(\alpha)_{\overline{i}} = \alpha_i$ for any $i \in Q_0$. Notice that if $V$ is an invariant representation of $Q$, then $\dim V_i = \dim V_{\sigma^{-1}(i)}$ for all $i \in Q_0$. As such, $\dim V_i = \dim V_j$ for all $i$ and $j$ in the same orbit. Thus, $\dimbar V \in (\mathbb{Z}^{Q_0})^\sigma$\index{$\dimbar$}\index{Graded dimension}. We have the following result due to Hubery.

\begin{theo}\textbf{\emph{\cite[Theorem 1]{hub}}}\label{theo:hubery}
  Let $(Q,\sigma)$ be a quiver with automorphism, $\Omega$ a valued quiver such that $(Q,\sigma)$ folds into $\Omega$, and $K$ an algebraically closed field of characteristic not dividing the order of $\sigma$.
  \begin{enumerate}
    \item The images under $\f$ of the graded dimensions of the invariant-indecomposable representations\index{Representation!invariant-indecomposable} of $Q$ are the positive roots\index{Root!positive} of $\Phi(\Omega)$.

    \item If $\f(\alpha)$ is a real positive root, then there is a unique invariant-indecomposable representation of $Q$ with graded dimension $\alpha$ (up to isomorphism). 
  \end{enumerate}
\end{theo}

Theorem \ref{theo:hubery} tells us that if the indecomposables of $\Q$ are determined by the positive roots of $\Phi(\Omega)$ (such as in the case of species of Dynkin or extended Dynkin type or $K$-species over finite fields), then finding all the indecomposables of $\Q$ reduces to finding the indecomposables of $Q$, which, in general, is an easier task.

One may wonder if there is a subcategory of $\R_K(Q)$\index{$\R_K(Q)$}, say $\R_K^\sigma(Q)$, whose objects are the invariant representations of $Q$, that is equivalent to $\R(\Q)$\index{$\R(\Q)$}. One needs to determine what the morphisms of this category should be. The most obvious choice is to let $\R_K^\sigma(Q)$ be the full subcategory of $\R_K(Q)$ whose objects are the invariant representations. This, however, does not work. The category $\R(\Q)$ is an abelian category (this follows from Proposition \ref{prop:equ}), but $\R_K^\sigma(Q)$, as we have defined it, is not. As the following example demonstrates, this category does not, in general, have kernels.

\begin{egg}
  Let $(Q,\sigma)$ be the following quiver with automorphism (where the dotted arrows represent the action of $\sigma$).

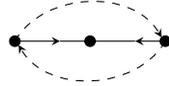
  \begin{center}
    \begin{tikzpicture}[>=stealth]
      \filldraw (0,0) circle (2pt) (1,0) circle (2pt) (2,0) circle (2pt);

      \draw[->] (0,0) -- (0.6,0); \draw (0.6,0) -- (1,0);
      \draw[->] (2,0) -- (1.6,0); \draw (1.6,0) -- (1,0);

      \draw[->, dashed, out=60, in=120] (0,0) to (1.95,0.05);
      \draw[->, dashed, out=-120, in=-60] (2,0) to (0.05,-0.05);
    \end{tikzpicture}
  \end{center}
  Let $V$ be the following invariant-indecomposable representation of $Q$.

    \begin{center}
    \begin{tikzpicture}[>=stealth]
      \draw (0,0) node {$K$} (1,0) node {$0$} (2,0) node {$K$};

      \draw[->] (0.3,0) -- (0.7,0);
      \draw[->] (1.7,0) -- (1.3,0);
    \end{tikzpicture}
  \end{center}
  Let $\varphi: V \rightarrow V$ be the morphism defined by

  \begin{center}
    \begin{tikzpicture}
      \draw (-1.5,0) node {$V$} (-1.5,-1.5) node {$V$} (-1.8,-0.75) node {$\varphi$} (0,0) node {$K$} (1.5,0) node {$0$} (3,0) node {$K$} (0,-1.5) node {$K$} (1.5,-1.5) node {$0$} (3,-1.5) node {$K$} (-0.3,-0.75) node {$0$} (1.2,-0.75) node {$0$} (2.7,-0.75) node {$1$};

      \draw[->] (-1.5,-0.3) -- (-1.5,-1.2); \draw[->] (0.3,0) -- (1.2,0); \draw[->] (2.7,0) -- (1.8,0); \draw[->] (0.3,-1.5) -- (1.2,-1.5); \draw[->] (2.7,-1.5) -- (1.8,-1.5); \draw[->] (0,-0.3) -- (0,-1.2); \draw[->] (1.5,-0.3) -- (1.5,-1.2); \draw[->] (3,-0.3) -- (3,-1.2);
    \end{tikzpicture}.
  \end{center}
  Then $\varphi$ is a morphism of representations since each of the squares in the diagram commutes. However, by a straightforward exercise in category theory, one can show that $\varphi$ does not have a kernel (in the category $\R_K^\sigma(Q)$).
\end{egg}

Therefore, if we define $\R_K^\sigma(Q)$ as a full subcategory of $\R_K(Q)$, it is not equivalent to $\R(\Q)$. It is possible that one could cleverly define the morphisms of $\R_K^\sigma(Q)$ to avoid this problem, however there are other obstacles to overcome. If $\R_K^\sigma(Q)$ and $\R(\Q)$ were equivalent, then there should be a bijective correspondence between the (isomorphism classes of) indecomposables in each category. Using the idea of folding, an invariant representation of $Q$ with graded dimension\index{Graded dimension} $\alpha$ should be mapped to a representation of $\Q$ with graded dimension $\f(\alpha)$\index{$\f$}. The following example illustrates the problem with this idea. Note that this example is similar to the example following Proposition 15 in \cite{hub}, however we approach it in a different fashion.

\begin{egg}
  Let $(Q,\sigma)$ be the following quiver with automorphism (again, the dotted arrows represent the action of $\sigma$).

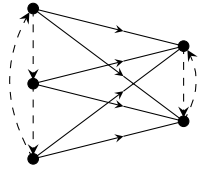
  \begin{center}
    \begin{tikzpicture}[>=stealth]
      \filldraw (0,0) circle (2pt) (0,1) circle (2pt) (0,2) circle (2pt) (2,0.5) circle (2pt) (2,1.5) circle (2pt);

      \draw[->] (0,0) -- (1.2,0.3); \draw[->] (0,0) -- (1.2,0.9); \draw[->] (0,1) -- (1.2,0.7); \draw[->] (0,1) -- (1.2,1.3); \draw[->] (0,2) -- (1.2,1.1); \draw[->] (0,2) -- (1.2,1.7);
      \draw (1.2,0.3) -- (2,0.5) (1.2,0.9) -- (2,1.5) (1.2,0.7) -- (2,0.5) (1.2,1.3) -- (2,1.5) (1.2,1.1) -- (2,0.5) (1.2,1.7) -- (2,1.5);

         \draw[decoration={markings,mark=at position 1 with {\arrow[thick]{>}}},postaction={decorate}, dashed] (0,2) -- (0,1.05); \draw[decoration={markings,mark=at position 1 with {\arrow[thick]{>}}},postaction={decorate}, dashed] (0,1) -- (0,0.05); \draw[decoration={markings,mark=at position 1 with {\arrow[thick]{>}}},postaction={decorate}, dashed, out=120, in=-120] (0,0) to (-0.05,1.95);
      \draw[decoration={markings,mark=at position 1 with {\arrow[thick]{>}}},postaction={decorate}, dashed] (2,1.5) -- (2,0.55); \draw[decoration={markings,mark=at position 1 with {\arrow[thick]{>}}},postaction={decorate}, dashed, out=60, in=-60] (2,0.5) to (2.05,1.45);
    \end{tikzpicture}
  \end{center}
  Then $(Q,\sigma)$ folds into the following absolute valued quiver.

  \begin{center}
    \begin{tikzpicture}[>=stealth]
      \draw (-1,0) node {$\Gamma$:};

      \filldraw (0,0) circle (2pt) (1,0) circle (2pt);

      \draw[->] (0,0) -- (0.6,0); \draw (0.6,0) -- (1,0);

      \draw (0,-0.3) node {$(3)$} (1,-0.3) node {$(2)$} (0.5,0.3) node {$(6)$};
    \end{tikzpicture}
  \end{center}
  One can easily check that $\beta=(1,1)$ is an imaginary root of $\Phi(\Gamma)$. The only $\alpha \in (\mathbb{Z}^{Q_0})^\sigma$ such that $\f(\alpha)=\beta$ is $\alpha=(1,1, \dots, 1)$. One can show using basic linear algebra that, while there are several non-isomorphic indecomposable representations of $Q$ with graded dimension $\alpha$ (after all, $\alpha$ is an imaginary root of $\Phi(Q)$), all such invariant representations are isomorphc to
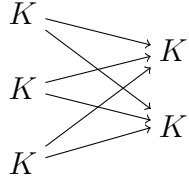
    \begin{center}
    \begin{tikzpicture}
      \filldraw (0,0) node {$K$} (0,1) node {$K$} (0,2) node {$K$} (2,0.5) node {$K$} (2,1.5) node {$K$};

      \draw[->] (0.3,0.075) -- (1.7,0.475); \draw[->] (0.3,0.225) -- (1.7,1.275); \draw[->] (0.3,0.925) -- (1.7,0.575); \draw[->] (0.3,1.075) -- (1.7,1.425); \draw[->] (0.3,1.775) -- (1.7,0.725); \draw[->] (0.3,1.925) -- (1.7,1.575);
    \end{tikzpicture}
  \end{center}
  where every arrow represents the identity map $\id_K$. Thus, we have a single isomorphism class of invariant-indecomposables with graded dimension $\alpha$.

  Now, construct a species of $\Gamma$. Let $\gamma=2^{1/6}$ and let $\Q$ be the $\mathbb{Q}$-species given by $\mathbb{Q}(\gamma^2) \xrightarrow{\mathbb{Q}(\gamma)} \mathbb{Q}(\gamma^3)$. Thus the underlying valued quiver of $\Q$ is $\Gamma$. There exists an indecomposable representation of $\Q$ with graded dimension $\beta$ -- in fact, there exists more than one.

  Let $V_1$ be the representation $\mathbb{Q}(\gamma^2) \xrightarrow{f_1} \mathbb{Q}(\gamma^3)$ where $f_1: \mathbb{Q}(\gamma) \otimes_{\mathbb{Q}(\gamma^2)} \mathbb{Q}(\gamma^2) \cong \mathbb{Q}(\gamma) \rightarrow \mathbb{Q}(\gamma^3)$ is the $\mathbb{Q}(\gamma^3)$-linear map defined by $1 \mapsto 1$, $\gamma \mapsto 0$ and $\gamma^2 \mapsto 0$.

  Let $V_2$ be the representation $\mathbb{Q}(\gamma^2) \xrightarrow{f_2} \mathbb{Q}(\gamma^3)$ where $f_2: \mathbb{Q}(\gamma) \otimes_{\mathbb{Q}(\gamma^2)} \mathbb{Q}(\gamma^2) \cong \mathbb{Q}(\gamma) \rightarrow \mathbb{Q}(\gamma^3)$ is the $\mathbb{Q}(\gamma^3)$-linear map defined by $1 \mapsto 1$, $\gamma \mapsto 1$ and $\gamma^2 \mapsto 0$.

  It is clear that $V_1$ and $V_2$ are indecomposable and $\dimbar V_1 = \dimbar V_2 = \beta$. One can also show that they are not isomorphic as representations of $\Q$. Hence, there are at least two isomorphism classes of indecomposable representations of $\Q$ with graded dimension $\beta$.

  Therefore, any functor $\R^\sigma_K(Q) \rightarrow \R(\Q)$ mapping invariant representations with graded dimension $\alpha$ to representations with graded dimension $\f(\alpha)$ cannot be essentially surjective, and thus cannot be an equivalence of categories.
\end{egg}

While the above example is not enough to conclude that the categories $\R^\sigma_K(Q)$ and $\R(\Q)$ are not equivalent, it is enough to deduce that one cannot obtain an equivalence via folding.

\section{Ringel-Hall Algebras}
In this section we define the Ringel-Hall algebra of a species (or quiver). We will construct the generic composition algebra of a species, which is obtained from a subalgebra of the Ringel-Hall algebra, and see that it is isomorphic to the positive part of the quantized enveloping algebra of the corresponding Kac-Moody Lie algebra (see Theorem \ref{theo:compalg}). We then give a similar interpretation of the whole Ringel-Hall algebra (see Theorem \ref{theo:fullrhalg}). For further details, see the expository paper by Schiffmann, \cite{schiff}.

We continue our assumption that all quivers/species have no oriented cycles. Also, we have seen in the last section (Proposition \ref{prop:equ}) that $\R(\Q)$ is equivalent to $T(\Q)$-mod, and so we will simply identify representations of $\Q$ with modules of $T(\Q)$.

\begin{defn}[Ringel-Hall algebra]\label{defn:rhalg}
  Let $\Q$ be an $\mathbb{F}_q$-species\index{Species!$K$-species}. Let $v=q^{1/2}$ and let $\A$ be an integral domain containing $\mathbb{Z}$ and $v,v^{-1}$. The \emph{Ringel-Hall algebra}\index{Ringel-Hall algebra}, which we will denote $\HH(\Q)$\index{$\HH(\Q)$}, is the free $\A$-module with basis the set of all isomorphism classes of finite-dimensional representations of $\Q$. Multiplication is given by\index{$\dimbar$}
  \[ [A][B]=v^{\langle \dimbar A,\dimbar B \rangle} \sum_{[C]} g_{AB}^C [C], \]
  where $g_{AB}^C$ is the number of subrepresentations (submodules) $X$ of $C$ such that $C/X \cong A$ and $X \cong B$ (as representations/modules) and $\langle-,-\rangle$ is the Euler form\index{Euler form} (see Definition \ref{defn:euler}).
\end{defn}

\begin{rmk}
  It is well-known (see, for example, \cite[Lemma 2.2]{r}) that
  \[ \langle \dimbar A, \dimbar B \rangle = \dim_{\mathbb{F}_q} \Hom_{T(\Q)}(A,B) - \dim_{\mathbb{F}_q} \Ext_{T(\Q)}(A,B). \]
  In many texts (for example \cite{green} or \cite{sv}) this is the way the form $\langle -,- \rangle$ is defined. Also, there does not appear to be a single agreed-upon name for this algebra; depending on the text, it may be called the \emph{twisted Hall algebra}, the \emph{Ringel algebra}, the \emph{twisted Ringel-Hall algebra}, etc. Regardless of the name one prefers, it is important not to confuse this algebra with the (untwisted) Hall algebra whose multiplication is given by $[A][B]=\sum_{[C]} g_{AB}^C [C]$.
\end{rmk}

\begin{defn}[Composition algebra]
  Let $\Q$ be an $\mathbb{F}_q$-species with underlying absolute valued quiver $\Gamma$ and $\mathbb{F}_q$-modulation $(K_i, M_\rho)_{i \in \Gamma_0, \rho \in \Gamma_1}$. The \emph{composition algebra}\index{Composition algebra}, $C=C(\Q)$, of $\Q$ is the $\A$-subalgebra of $\HH(\Q)$ generated by the isomorphism classes of the simple representations of $\Q$. Since we assume $\Gamma$ has no oriented cycles, this means $C$ is generated by the $[S_i]$ for $i \in \Gamma_0$ where $S_i=((S_i)_j,(S_i)_\rho)_{j \in \Gamma_0, \rho \in \Gamma_1}$ is given by
  \[ (S_i)_j = \begin{cases} K_i, & \text{if } i=j, \\ 0, & \text{if } i\ne j, \end{cases} \qquad \text{and} \qquad (S_i)_\rho = 0 \quad \text{for all } \rho \in \Gamma_1. \]
\end{defn}

Let $\mathcal{S}$ be a set of finite fields $K$ such that $\{ |K| \; | \; K \in \mathcal{S} \}$ is infinite. Let $v_K = |K|^{1/2}$ for each $K \in \mathcal{S}$. Write $C_K$ for the composition algebra of $\Q$ for each finite field $K$ in $\mathcal{S}$ and $[S_i^{(K)}]$ for the corresponding generators. Let $C$ be the subring of $\Pi_{K \in \mathcal{S}} C_K$ generated by $\mathbb{Q}$ and the elements
\[ t = (t_K)_{K \in \mathcal{S}}, \; t_K = v_K, \]
\[ t^{-1} = (t^{-1}_K)_{K \in \mathcal{S}}, \; t^{-1}_K = v^{-1}_K, \]
\[ u_i = (u_i^{(K)})_{K \in \mathcal{S}}, \; u_i^{(K)} = [S_i^{(K)}]. \]
So $t$ lies in the centre of $C$ and, because there are infinitely many $v_K$, $t$ does not satisfy $p(t)=0$ for any nonzero polynomial $p(T)$ in $\mathbb{Q}[T]$. Thus, we may view $C$ as the $\A$-algebara generated by the $u_i$, where $\A = \mathbb{Q}[t,t^{-1}]$ with $t$ viewed as an indeterminate.

\begin{defn}[Generic composition algebra]
  Using the notation above, the $\mathbb{Q}(t)$-algebra $C^* = \mathbb{Q}(t) \otimes_\A C$ is called the \emph{generic composition algebra}\index{Composition algebra!generic} of $\Q$. We write $u^*_i = 1 \otimes u_i$ for $i \in \Gamma_0$.
\end{defn}

Let $\Q$ be an $\mathbb{F}_q$-species with underlying absolute valued quiver $\Gamma$. Let $(c_{ij})$ be the generalized Cartan matrix associated to $\Gamma$ and let $\g$ be its associated Kac-Moody Lie algebra\index{Kac-Moody Lie algebra} (recall Definitions \ref{defn:cartan} and \ref{defn:kac-moody}). Let $U_t(\g)$ be the quantized enveloping algebra\index{Quantized enveloping algebra} of $\g$ and let $U_t(\g) = U_t^+(\g) \otimes U_t^0(\g) \otimes U_t^-(\g)$ be its triangular decomposition (see \cite[Chapter 3]{lusz}). We call $U_t^+(\g)$ the \emph{positive part} of $U_t(\g)$; it is the $\mathbb{Q}(t)$-algebra generated by elements $E_i$, $i \in \Gamma_0$, modulo the quantum Serre relations
\[ \sum_{p=0}^{1-c_{ij}} (-1)^p \begin{bmatrix} 1-c_{ij} \\ p \end{bmatrix} E_i^p E_j E_i^{1-c_{ij}-p} \qquad \text{for all } i \ne j, \]
where
\[ \begin{bmatrix} m \\ p \end{bmatrix} = \frac{[m]!}{[p]! [m-p]!},\]
\[ [n]=\frac{t^n - t^{-n}}{t - t^{-1}}, \qquad [n]! = [1] [2] \cdots [n].\]
In \cite{green}, Green was able to show that $C^*$ and $U_t^+(\g)$ are canonically isomorphic. Of course, in his paper, Green speaks of modules of hereditary algebras\index{Hereditary algebra} over a finite field $K$ rather than representations of $K$-species, but as we have seen, these two notions are equivalent.

\begin{theo}\cite[Theorem 3]{green}\label{theo:compalg}
  Let $\Q$ be an $\mathbb{F}_q$-species with underlying absolute valued quiver $\Gamma$ and let $\g$ be its associated Kac-Moody Lie algebra. Then, there exists a $\mathbb{Q}(t)$-algebra isomorphism $U_t^+(\g) \rightarrow C^*$ which takes $E_i \mapsto u^*_i$ for all $i \in \Gamma_0$.
\end{theo}

\begin{rmk}
  This result by Green is actually a generalization of an earlier result by Ringel in \cite[p. 400]{r3} and \cite[Theorem 7]{r4} who proved Theorem 7.5 in the case that $\Q$ is of finite representation type.
\end{rmk}

Theorem \ref{theo:compalg} gives us an interpretation of the composition algebra in terms of the quantized enveloping algebra of the corresponding Kac-Moody Lie algebra. Later, Sevenhant and Van Den Bergh were able to give a similar interpretation of the whole Ringel-Hall algebra. For this, however, we need the concept of a \emph{generalized Kac-Moody Lie algebra}, which was first defined by Borcherds in \cite{bor}. Though some authors have used slightly modified definitions of generalized Kac-Moody Lie algebras over the years, we use here Borcherds's original definition (in accordance with Sevenhant and Van Den Bergh in \cite{sv}).

\begin{defn}[Generalized Kac-Moody Lie algebra]\label{defn:genkac-moody}
  Let $H$ be a real vector space with symmetric bilinear product $(-,-): H \times H \rightarrow \mathbb{R}$. Let $\mathcal{I}$ be a countable (but possibly infinite) set and $\{ h_i \}_{i \in \mathcal{I}}$ be a subset of $H$ such that $(h_i,h_j) \le 0$ for all $i \ne j$ and $c_{ij}=2(h_i,h_j)/(h_i,h_i)$ is an integer if $(h_i,h_i) > 0$. Then, the \emph{generalized Kac-Moody Lie algebra}\index{Kac-Moody Lie algebra!generalized} associated to $H$, $\{ h_i \}_{i \in \mathcal{I}}$ and $(-,-)$ is the Lie algebra (over a field of characteristic $0$ containing an isomorphic copy of $\mathbb{R}$) generated by $H$ and elements $e_i$ and $f_i$ for $i \in \mathcal{I}$ whose product is defined by:

  \begin{itemize}
    \item $[h,h'] = 0$ for all $h$ and $h'$ in $H$,
    \item $[h,e_i] = (h,h_i)e_i$ and $[h,f_i] = -(h,h_i)f_i$ for all $h \in H$ and $i \in \mathcal{I}$,
    \item $[e_i,f_i] = h_i$ for each $i$ and $[e_i,f_j]=0$ for all $i \ne j$,
    \item if $(h_i,h_i) > 0$, then $(\ad e_i)^{1-c_{ij}} (e_j) = 0$ and $(\ad f_i)^{1-c_{ij}}(f_j) = 0$ for all $i \ne j$,
    \item if $(h_i,h_i)=0$, then $[e_i,e_j]=[f_i,f_j]=0$.
  \end{itemize}
\end{defn}

\begin{rmk}
  Generalized Kac-Moody Lie algebras are similar to Kac-Moody Lie algebras. The main difference is that generalized Kac-Moody Lie algebras (may) contain simple imaginary roots\index{Root!imaginary} (corresponding to the $h_i$ with $(h_i,h_i) \le 0$).
\end{rmk}

Let $\g$ be a generalized Kac-Moody Lie algebra (with the notation of Definition \ref{defn:genkac-moody}) and let $v \ne 0$ be an element of the base field such that $v$ is not a root of unity. Write $d_i = (h_i,h_i)/2$ for the $i \in \mathcal{I}$ such that $(h_i,h_i)>0$. The \emph{quantized enveloping algebra} $U_v(\g)$ can be defined in the same way as the quantized enveloping algebra\index{Quantized enveloping algebra} of a Kac-Moody Lie algebra (see Section 2 of \cite{sv}). The positive part of $U_v^+(\g)$ is the $\A$-algebra generated by elements $E_i$, $i \in \mathcal{I}$, modulo the quantum Serre relations
\[ \sum_{p=0}^{1-c_{ij}} (-1)^p \begin{bmatrix} 1-c_{ij} \\ p \end{bmatrix}_{d_i} E_i^p E_j E_i^{1-c_{ij}-p} \qquad \text{for all } i \ne j, \text{ with } (h_i,h_i)>0 \]
and
\[ E_i E_j - E_j E_i \qquad \text{if } (h_i,h_j) = 0,\]
where
\[ \begin{bmatrix} m \\ p \end{bmatrix}_{d_i} = \frac{[m]_{d_i}!}{[p]_{d_i}! [m-p]_{d_i}!},\]
\[ [n]_{d_i}=\frac{(v^{d_i})^n - (v^{d_i})^{-n}}{v^{d_i} - v^{-d_i}}, \qquad [n]_{d_i}! = [1]_{d_i} [2]_{d_i} \cdots [n]_{d_i}. \]

Let $n$ be a positive integer and let $\{ e_i \}_{i=1}^n$ be the standard basis of $\mathbb{Z}^n$. Let $v \in \mathbb{R}$ such that $v > 1$ and $\A$ be as in Definition \ref{defn:rhalg}. Suppose we have the following:

\begin{enumerate}
  \item An $\mathbb{N}^n$-graded $\A$-algebra $A$ such that:
    \begin{enumerate}
      \item $A_0=\A$,
      \item $\dim_\A A_\alpha < \infty$ for all $\alpha \in \mathbb{N}^n$,
      \item $A_{e_i} \ne 0$ for all $1\le i \le n$.
    \end{enumerate}

  \item A symmetric positive definite bilinear form $[-,-]:A \times A \rightarrow \A$ such that $[A_\alpha,A_\beta]=0$ if $\alpha \ne \beta$ and $[1,1]=1$ (here we assume $[a,a] \in \mathbb{R}$ for all $a \in A)$.

  \item A symmetric bilinear form $(-,-): \mathbb{R}^n \times \mathbb{R}^n \rightarrow \mathbb{R}$ such that $(e_i,e_i) > 0$ for all $1 \le i \le n$ and $c_{ij}=2(e_i,e_j)/(e_i,e_i)$ is a generalized Cartan matrix\index{generalized Cartan matrix} as in Definition 6.9.

  \item The tensor product $A \otimes_\A A$ can be made into an algebra via the rule
  \[ (a \otimes b)(c \otimes d) = v^{(\deg(b),\deg(c))}(ac \otimes bd), \]
  for homogeneous $a,b,c,d$. (here $\deg(x)=\alpha$ if $x \in A_\alpha$). We assume that there is an $\A$-algebra homomorphism $\boldsymbol{\delta} : A \rightarrow A \otimes_\A A$ which is adjoint under $[-,-]$ to the multiplication (that is, $[\boldsymbol{\delta}(a),b \otimes c]_{A \otimes A}=[a,bc]_{A}$ where $[a \otimes b, c \otimes d]_{A \otimes A} = [a,c]_A [b,d]_A$).
\end{enumerate}

\begin{prop}\textbf{\emph{\cite[Proposition 3.2]{sv}}}
  Under the above conditions, $A$ is isomorphic (as an algebra) to the positive part of the quantized enveloping algebra of a generalized Kac-Moody Lie algebra.
\end{prop}

The Ringel-Hall algebra\index{Ringel-Hall algebra}, $\HH(\mathcal{Q})$\index{$\HH(\Q)$}, of a species $\mathcal{Q}$ is $\mathbb{N}^n$-graded by associating to each representation its graded dimension, hence $\HH(\Q)$ satisfies Condition 1 above. Moreover, the symmetric Euler form satisfies Condition 3 (if we extend it to all $\mathbb{R}^{\Gamma_0}$). Following Green in \cite{green}, we define
\[ \boldsymbol{\delta}\left([A]\right)=\sum_{[B],[C]} v^{\langle \dimbar B,\dimbar C \rangle} g_{BC}^A \dfrac{|\Aut(B)||\Aut(C)|}{|\Aut(A)|} \left([B] \otimes [C]\right) \]
and
\[([A],[B])_{\HH(\Q)} = \dfrac{\delta_{[A],[B]}}{|\Aut(A)|}. \]
In \cite[Theorem 1]{green}, Green shows that $(-,-)_{\HH(\Q)}$ satisfies Condition 2 and that $\boldsymbol{\delta}$ satisfies Condition 4. Hence, we have the following.

\begin{theo}\textbf{\emph{\cite[Theorem 1.1]{sv}}}\label{theo:fullrhalg}
  Let $\Q$ be an $\mathbb{F}_q$-species. Then, $\HH(\Q)$ is the positive part of the quantized enveloping algebra\index{Quantized enveloping algebra} of a generalized Kac-Moody algebra\index{Kac-Moody Lie algebra!generalized}.
\end{theo}

\begin{rmk}
  In their paper, Sevenhant and Van Den Bergh state Theorem \ref{theo:fullrhalg} only for the Ringel-Hall algebra of a quiver, but none of their arguments depend on having a quiver rather than a species. Indeed, many of their arguments are based on those of Green in \cite{green}, which are valid for hereditary algebras\index{Hereditary algebra}. Moreover, Sevenhant and Van Den Bergh define the Ringel-Hall algebra to be an algebra opposite to the one we defined in Definition \ref{defn:rhalg} (our definition, which is the one used by Green, seems to be the more standard definition). This does not affect any of the arguments presented. 
\end{rmk}

\bibliographystyle{plain}

\bibliography{article}		

\begin{thebibliography}{10}

\bibitem{ass}
Ibrahim Assem, Daniel Simson, and Andrzej Skowro{\'n}ski.
\newblock {\em Elements of the representation theory of associative algebras.
  {V}ol. 1}, volume~65 of {\em London Mathematical Society Student Texts}.
\newblock Cambridge University Press, Cambridge, 2006.
\newblock Techniques of representation theory.

\bibitem{aus}
Maurice Auslander.
\newblock On the dimension of modules and algebras. {III}. {G}lobal dimension.
\newblock {\em Nagoya Math. J.}, 9:67--77, 1955.

\bibitem{b}
D.~J. Benson.
\newblock {\em Representations and cohomology. {I}}, volume~30 of {\em
  Cambridge Studies in Advanced Mathematics}.
\newblock Cambridge University Press, Cambridge, second edition, 1998.
\newblock Basic representation theory of finite groups and associative
  algebras.

\bibitem{bor}
Richard Borcherds.
\newblock Generalized {K}ac-{M}oody algebras.
\newblock {\em J. Algebra}, 115(2):501--512, 1988.

\bibitem{bourbaki}
N.~Bourbaki.
\newblock {\em {\'E}l{\'e}ments de math{\'e}matique: {A}lg{\`e}bre, Chapitre 4
  {\`a} 7}.
\newblock Springer, 2007.

\bibitem{brauer}
Richard Brauer.
\newblock {\em Collected papers. {V}ol. {I}}, volume~17 of {\em Mathematicians
  of Our Time}.
\newblock MIT Press, Cambridge, Mass., 1980.
\newblock Theory of algebras, and finite groups, Edited and with an
  introduction by Paul Fong and Warren J. Wong, With a biography by J. A.
  Green, With an introduction by O. Goldman.

\bibitem{cnp}
Xueqing Chen, Ki-Bong Nam, and Tom{\'a}\v{s} Posp{\'i}chal.
\newblock Quivers and representations.
\newblock In {\em Handbook of algebra. {V}ol. 6}, volume~6 of {\em Handb.
  Algebr.}, pages 507--561. Elsevier/North-Holland, Amsterdam, 2009.

\bibitem{dx1}
B.~Deng and J.~Xiao.
\newblock The {R}ingel-{H}all interpretation to a conjecture of {K}ac.
\newblock preprint.

\bibitem{ddpw}
Bangming Deng, Jie Du, Brian Parshall, and Jianpan Wang.
\newblock {\em Finite dimensional algebras and quantum groups}, volume 150 of
  {\em Mathematical Surveys and Monographs}.
\newblock American Mathematical Society, Providence, RI, 2008.

\bibitem{dr2}
Vlastimil Dlab and Claus~Michael Ringel.
\newblock On algebras of finite representation type.
\newblock {\em J. Algebra}, 33:306--394, 1975.

\bibitem{dr}
Vlastimil Dlab and Claus~Michael Ringel.
\newblock Indecomposable representations of graphs and algebras.
\newblock {\em Mem. Amer. Math. Soc.}, 6(173):v+57, 1976.

\bibitem{dk}
Yurij~A. Drozd and Vladimir~V. Kirichenko.
\newblock {\em Finite-dimensional algebras}.
\newblock Springer-Verlag, Berlin, 1994.
\newblock Translated from the 1980 Russian original and with an appendix by
  Vlastimil Dlab.

\bibitem{g}
Peter Gabriel.
\newblock Indecomposable representations. {II}.
\newblock In {\em Symposia {M}athematica, {V}ol. {XI} ({C}onvegno di {A}lgebra
  {C}ommutativa, {INDAM}, {R}ome, 1971)}, pages 81--104. Academic Press,
  London, 1973.

\bibitem{green}
James~A. Green.
\newblock Hall algebras, hereditary algebras and quantum groups.
\newblock {\em Invent. Math.}, 120(2):361--377, 1995.

\bibitem{hub}
Andrew Hubery.
\newblock Quiver representations respecting a quiver automorphism: a
  generalisation of a theorem of {K}ac.
\newblock {\em J. London Math. Soc. (2)}, 69(1):79--96, 2004.

\bibitem{jacobson}
Nathan Jacobson.
\newblock {\em Basic algebra. {II}}.
\newblock W. H. Freeman and Company, New York, second edition, 1989.

\bibitem{jans}
James~P. Jans.
\newblock On the indecomposable representations of algebras.
\newblock {\em Ann. of Math. (2)}, 66:418--429, 1957.

\bibitem{kac2}
V.~G. Kac.
\newblock Infinite root systems, representations of graphs and invariant
  theory.
\newblock {\em Invent. Math.}, 56(1):57--92, 1980.

\bibitem{tyl}
T.~Y. Lam.
\newblock {\em A first course in noncommutative rings}, volume 131 of {\em
  Graduate Texts in Mathematics}.
\newblock Springer-Verlag, New York, 1991.

\bibitem{me}
Joel Lemay.
\newblock Valued graphs and the representation theory of {L}ie algebras.
\newblock Master's thesis, University of Ottawa, 2011.
\newblock http://www.ruor.uottawa.ca/en/handle/10393/20168.

\bibitem{lusz}
George Lusztig.
\newblock {\em Introduction to quantum groups}.
\newblock Modern Birkh{\"a}user Classics. Birkh{\"a}user/Springer, New York,
  2010.
\newblock Reprint of the 1994 edition.

\bibitem{bn}
Bertrand Nguefack.
\newblock A non-simply laced version for cluster structures on 2-{C}alabi-{Y}au
  categories.
\newblock {\em Preprint}, 2009.

\bibitem{r2}
Claus~Michael Ringel.
\newblock Representations of {$K$}-species and bimodules.
\newblock {\em J. Algebra}, 41(2):269--302, 1976.

\bibitem{r}
Claus~Michael Ringel.
\newblock Hall algebras and quantum groups.
\newblock {\em Invent. Math.}, 101(3):583--591, 1990.

\bibitem{r3}
Claus~Michael Ringel.
\newblock From representations of quivers via {H}all and {L}oewy algebras to
  quantum groups.
\newblock In {\em Proceedings of the {I}nternational {C}onference on {A}lgebra,
  {P}art 2 ({N}ovosibirsk, 1989)}, volume 131 of {\em Contemp. Math.}, pages
  381--401, Providence, RI, 1992. Amer. Math. Soc.

\bibitem{r4}
Claus~Michael Ringel.
\newblock The {H}all algebra approach to quantum groups.
\newblock In {\em X{I} {L}atin {A}merican {S}chool of {M}athematics ({S}panish)
  ({M}exico {C}ity, 1993)}, volume~15 of {\em Aportaciones Mat. Comun.}, pages
  85--114. Soc. Mat. Mexicana, M{\'e}xico, 1995.

\bibitem{schiff}
O.~Schiffmann.
\newblock Lectures on {H}all algebras.
\newblock arXiv:math/0611617v2.

\bibitem{sv}
Bert Sevenhant and Michel {Van Den Bergh}.
\newblock A relation between a conjecture of {K}ac and the structure of the
  {H}all algebra.
\newblock {\em J. Pure Appl. Algebra}, 160(2-3):319--332, 2001.

\bibitem{weil}
Andr{\'e} Weil.
\newblock {\em Basic number theory}.
\newblock Springer-Verlag, New York, third edition, 1974.
\newblock Die Grundlehren der Mathematischen Wissenschaften, Band 144.

\bibitem{yoshii}
Tensho Yoshii.
\newblock On algebras of bounded representation type.
\newblock {\em Osaka Math. J.}, 8:51--105, 1956.

\end{thebibliography}

\end{document}